\documentclass{amsart}
\usepackage{graphicx}
\usepackage{amscd}
\usepackage{amsmath}%
\usepackage{amsfonts}%
\usepackage{amssymb}
\theoremstyle{plain}
\newtheorem{theorem}{Theorem}[section]

\newtheorem{lemma}[theorem]{Lemma}

\newtheorem{proposition}[theorem]{Proposition}
\theoremstyle{definition}
\newtheorem{example}[theorem]{Example}

\newtheorem{remark}[theorem]{Remark}

\numberwithin{equation}{section}

\begin{document}
\title[Solution of the Truncated Hyperbolic Moment Problem]{Solution of the Truncated Hyperbolic \\Moment Problem}
\author{Ra\'{u}l E. Curto}
\address{Department of Mathematics\newline The University of Iowa\newline Iowa City, IA
52242-1419\newline USA}
\email{rcurto@math.uiowa.edu}
\author{Lawrence A. Fialkow}
\address{Department of Computer Science\newline State University of New York\newline
New Paltz, NY 12561\newline USA}
\email{fialkowl@newpaltz.edu}
\thanks{The first-named author's research was partially supported by NSF Research
Grants DMS-9800931 and DMS-0099357. \ The second-named author's research was
partially supported by NSF Research Grant DMS-0201430. \ The second-named
author was also partially supported by the State University of New York at New
Paltz Research and Creative Projects Award Program.}

\begin{abstract}
Let $Q(x,y)=0$ be an hyperbola in the plane. \ Given real numbers $\beta
\equiv\beta^{\left(  2n\right)  }=\{\beta_{ij}\}_{i,j\geq0,i+j\leq2n}$, with
$\beta_{00}>0$, the \emph{truncated }$Q$\textit{-}\emph{hyperbolic moment
problem} for $\beta$ entails finding necessary and sufficient conditions for
the existence of a positive Borel measure $\mu$, supported in $Q(x,y)=0$, such
that $\beta_{ij}=\int y^{i}x^{j}\,d\mu\;\;(0\leq i+j\leq2n)$. We prove that
$\beta$ admits a $Q$\textit{-representing measure} $\mu$ (as above) if and
only if the associated moment matrix $\mathcal{M}(n)(\beta)$ is positive
semidefinite, recursively generated, has a column relation $Q(X,Y)=0$, and the
algebraic variety $\mathcal{V}(\beta)$ associated to $\beta$ satisfies
$\operatorname*{card}\mathcal{V}(\beta)\geq\operatorname*{rank}\mathcal{M}%
(n)(\beta)$. \ In this case, $\operatorname*{rank}\mathcal{M}(n)\leq2n+1$; if
$\operatorname*{rank}\mathcal{M}(n)\leq2n$, then $\beta$ admits a
$\operatorname*{rank}\mathcal{M}(n)$-atomic (minimal) $Q$-representing
measure; if $\operatorname*{rank}\mathcal{M}(n)=2n+1$, then $\beta$ admits a
$Q$-representing measure $\mu$ satisfying $2n+1\leq\operatorname*{card}%
\operatorname*{supp}\mu\leq2n+2$.

\end{abstract}
\maketitle

\section{\label{Int}Introduction}

Let $Q(x,y)=0$ denote an hyperbola in the plane. \ Given a real sequence
$\beta\equiv\beta^{(2n)}:\beta_{00},\beta_{01},\beta_{10},...,$ $\beta
_{0,2n},\beta_{2n,0}$, with $\beta_{00}>0$, we seek concrete necessary and
sufficient conditions so that there exists a positive Borel measure $\mu$ on
$\mathbb{R}^{2}$ satisfying
\begin{equation}
\beta_{ij}=\int y^{i}x^{j}\,d\mu\qquad(i,j\geq0,\;i+j\leq2n) \label{11}%
\end{equation}
and
\begin{equation}
\operatorname*{supp}\mu\subseteq\mathcal{Z}(Q):=\{(x,y)\in\mathbb{R}%
^{2}:Q(x,y)=0\}; \label{12}%
\end{equation}
a measure $\mu$ satisfying (\ref{11}) is a \textit{representing measure} for
$\beta$; $\mu$ is a $Q$\textit{-representing measure} if it satisfies
(\ref{11}) and (\ref{12}).

Our criterion for the existence of representing measures is expressed in terms
of algebraic and geometric properties of the \textit{moment matrix}
$\mathcal{M}(n)\equiv\mathcal{M}(n)(\beta)$ that we next describe. \ The size
of $\mathcal{M}(n)$ is $m(n):=(n+1)(n+2)/2$, with rows and columns indexed as
$1,X,Y,X^{2},YX,Y^{2},...,X^{n},YX^{n-1},...,Y^{n}$; the entry in row
$Y^{i}X^{j}$, column $Y^{k}X^{\ell}$ of $\mathcal{M}(n)$ is $\beta
_{i+k,j+\ell}$.

Let $\mathbb{R}_{n}[x,y]$ denote the space of real polynomials of degree at
most $n$ in two variables, and let $\mathcal{B}_{n}$ denote the basis
consisting of the monomials in degree-lexicographic order, i.e.,
$\mathcal{B}_{n}:1,x,y,x^{2},yx,y^{2},...,x^{n},yx^{n-1},...,y^{n}$. \ For
$p\in\mathbb{R}_{2n}[x,y]$, $p(x,y)\equiv\sum a_{rs}y^{r}x^{s}$, let $\hat
{p}:=(a_{rs})$ denote the coefficient vector of $p$ with respect to
$\mathcal{B}_{2n}$. \ Further, let $L_{\beta}:\mathbb{R}_{2n}[x,y]\rightarrow
\mathbb{R}$ be the \textit{Riesz functional} defined by $L_{\beta}(p):=\sum
a_{rs}\beta_{rs}$; then $\mathcal{M}(n)$ is uniquely determined by
\begin{equation}
\left(  \mathcal{M}(n)\hat{p},\hat{q}\right)  :=L_{\beta}(pq)\;\;(p,q\in
\mathbb{R}_{n}[x,y]). \label{13}%
\end{equation}
In particular, if $\mu$ is a representing measure for $\beta$, then $\left(
\mathcal{M}(n)\hat{p},\hat{p}\right)  =L_{\beta}(p^{2})=\int p^{2}\;d\mu\geq
0$. \ Since $\mathcal{M}(n)$ is real symmetric, it follows that $\left(
\mathcal{M}(n)(\hat{p}+i\hat{q}),\hat{p}+i\hat{q}\right)  \geq0$, whence
\begin{equation}
\mathcal{M}(n)\geq0 \label{14}%
\end{equation}
(i.e., $\mathcal{M}(n)$ is a positive semi-definite operator on $\mathbb{C}%
^{m(n)}$).

For $p\in\mathbb{R}_{n}[x,y]$, $p(x,y)\equiv\sum a_{ij}y^{i}x^{j}$, we define
an element $p(X,Y)$ of $\mathcal{C}_{\mathcal{M}(n)}$, the column space of
$\mathcal{M}(n)$, by $p(X,Y):=\sum a_{ij}Y^{i}X^{j}$; for polynomials $p$ and
$q$ with $\deg p+\deg q\leq n$, we also write $p(X,Y)q(X,Y)$ for $(pq)(X,Y)$.
\ It follows from \cite[Proposition 3.1]{tcmp1} that if $\mu$ is a
representing measure for $\beta$, then
\begin{equation}
\text{for }p\in\mathbb{R}_{n}[x,y],\;p(X,Y)=0\Leftrightarrow
\operatorname*{supp}\mu\subseteq\mathcal{Z}(p).\ \label{15}%
\end{equation}
It follows immediately from (\ref{15}) that if $\beta$ has a representing
measure, then $\mathcal{M}(n)$ is \textit{recursively generated }in the
following sense:
\begin{equation}
p,q,pq\in\mathbb{R}_{n}[x,y],\;p\left(  X,Y\right)  =0\Longrightarrow\left(
pq\right)  \left(  X,Y\right)  =0. \label{16}%
\end{equation}
We define the \textit{variety} of $\mathcal{M}(n)$ (or of $\beta$) by
$\mathcal{V}(\mathcal{M}(n)):=\bigcap_{\substack{p\in\mathbb{R}_{n}[x,y]
\\p(X,Y)=0}}\mathcal{Z}(p)$; \cite[Proposition 3.1 and Corollary 3.7]{tcmp1}
implies that if $\mu$ is a representing measure for $\beta$, then
$\operatorname*{supp}\mu\subseteq\mathcal{V}(\mathcal{M}(n))$ and
$\operatorname*{rank}\mathcal{M}(n)\leq\operatorname*{card}%
\operatorname*{supp}\mu\leq\operatorname*{card}\mathcal{V}(\mathcal{M}(n))$,
whence $\mathcal{M}(n)$ satisfies the \textit{variety condition}
\begin{equation}
\operatorname*{rank}\mathcal{M}(n)\leq\operatorname*{card}\mathcal{V}%
(\mathcal{M}(n)). \label{17}%
\end{equation}
In the sequel we repeatedly reply on the following basic result of
\cite[Theorem 5.13]{tcmp1}:
\begin{gather}
\beta\text{ admits a }\operatorname*{rank}\mathcal{M}(n)\text{-atomic
(minimal) representing measure }\label{exp18}\\
\text{if and only if }\mathcal{M}(n)\geq0\text{ and }\mathcal{M}(n)\text{
admits an extension to }\nonumber\\
\text{a (necessarily positive) moment matrix }\mathcal{M}(n+1)\text{
satisfying }\nonumber\\
\operatorname*{rank}\mathcal{M}(n+1)=\operatorname*{rank}\mathcal{M}%
(n);\nonumber
\end{gather}
we refer to such an extension as a \textit{flat extension}. \ 

Our main result shows that properties (\ref{14}) - (\ref{17}) completely
characterize the existence of $Q$-representing measures, as follows.

\begin{theorem}
\label{firstmain}Let $Q(x,y)=0$ be an hyperbola in the plane. A sequence
$\beta\equiv\beta^{(2n)}$ has a representing measure supported in $Q(x,y)=0 $
if and only if $\mathcal{M}(n)$ is positive semi-definite, recursively
generated, $Q(X,Y)=0$ in $\mathcal{C}_{\mathcal{M}(n)}$, and
$\operatorname*{rank}\mathcal{M}(n)\leq\operatorname*{card}\mathcal{V}%
(\mathcal{M}(n))$. \ In this case, $\operatorname*{rank}\mathcal{M}%
(n)\leq2n+1$; if $\operatorname*{rank}\mathcal{M}(n)\leq2n$, then there is a
$\operatorname*{rank}\mathcal{M}(n)$-atomic $Q$-representing measure, while if
$\operatorname*{rank}\mathcal{M}(n)=2n+1 $, there is a $Q$-representing
measure $\mu$ for which $2n+1\leq\operatorname*{card}\operatorname*{supp}%
\mu\leq2n+2$.
\end{theorem}

Consider the following property for a polynomial $P\in\mathbb{R}_{n}[x,y]$:
\begin{gather}
\beta\equiv\beta^{(2n)}\text{ has a representing measure supported in
}\mathcal{Z}(P)\text{ if and only if }\tag{$A_n^\prime$}\\
\mathcal{M}(n)(\beta)\text{ is positive semi-definite, recursively generated,
}\nonumber\\
P(X,Y)=0\text{ in }\mathcal{C}_{\mathcal{M}(n)}\text{, and }%
\operatorname*{rank}\mathcal{M}(n)\leq\operatorname*{card}\mathcal{V}%
(\mathcal{M}(n)).\nonumber
\end{gather}
Polynomials which satisfy $(A_{n}^{\prime})$ form an attractive class, because
if $P$ satisfies $(A_{n}^{\prime})$, then the degree-$2n$ moment problem on
$P(x,y)=0$ can be solved by concrete tests involving only elementary linear
algebra and the calculation of roots of polynomials. \ Theorem \ref{firstmain}
shows that each hyperbolic polynomial satisfies $(A_{n}^{\prime})$ for
$n\geq2$. \ Moreover, $P$ satisfies $(A_{n}^{\prime})$ for $n\geq\deg P$ if
$P$ represents a line \cite{tcmp2}, ellipse \cite{tcmp5}, or parabola
\cite{tcmp7}. \ These results together yield the following.

\begin{theorem}
\label{thm12}If $\deg P\leq2$, then $P$ satisfies $(A_{n}^{\prime})$ for every
$n\geq\deg P$.
\end{theorem}%

\noindent
Despite Theorem \ref{thm12}, there are differences between the parabolic and
elliptic moment problems and the hyperbolic problem. \ In the former cases,
the conditions of $(A_{n}^{\prime})$ always imply the existence of a
$\operatorname*{rank}\mathcal{M}(n)$-atomic representing measure,
corresponding to a flat extension of $\mathcal{M}(n)$; for this reason,
positive Borel measures supported on these curves always admit
\textit{Gaussian} cubature rules, i.e., $\operatorname*{rank}\mathcal{M}%
(n)$-atomic cubature rules of degree $2n$ (cf. \cite{FP}). \ By contrast, in
the hyperbolic case, minimal representing measures $\mu$ sometimes entail
$\operatorname*{card}\operatorname*{supp}\mu>\operatorname*{rank}%
\mathcal{M}(n)$ (and Gaussian cubature rules may fail to exist; cf. Example
\ref{exdegcaseIII}).

The preceding results are part of a general study of truncated moment problems
that we initiated in \cite{tcmp1}, and are also related to the classical full
moment problem, where moments of all orders are prescribed, i.e., $\beta
\equiv\beta^{(\infty)}=(\beta_{ij})_{i,j\geq0}$ (cf. \cite{AhKr}, \cite{Akh},
\cite{KrNu}, \cite{PuVa}, \cite{Sch}, \cite{ShTa}, \cite{StSz2}). \ Theorem
\ref{thm12} is motivated in part by results of J. Stochel \cite{Sto1}, who
solved the full moment problem on planar curves of degree at most $2$.
\ Paraphrasing \cite{Sto1} (i.e., translating from the language of
\textit{moment sequences} into the language of moment matrices), we consider
the following property of a polynomial $P$:
\begin{gather}
\beta^{(\infty)}\text{ has a representing measure supported in }%
P(x,y)=0\tag{$A$}\\
\text{if and only if }\mathcal{M}(\infty)(\beta)\geq0\text{ and }%
P(X,Y)=0\text{ in }\mathcal{C}_{\mathcal{M}(\infty)}.\nonumber
\end{gather}

\begin{theorem}
\label{thm13}(Stochel \cite{Sto1}) \ If $\deg P\leq2$, then $P$ satisfies
$(A)$.
\end{theorem}%

\noindent
In \cite{Sto1}, Stochel also proved that there exist polynomials of degree $3
$ that do not satisfy $(A)$; Stochel and F. Szafraniec \cite{StSz1} proved
that there are polynomials of arbitrarily large degree that satisfy $(A)$ (cf.
\cite{FiaNew}). \ Whether there exists a polynomial $P$ such that $P$ fails to
satisfy $(A_{n}^{\prime})$ for some $n\geq\deg P$ is an open problem (cf.
Section \ref{application}).

The link between the truncated and full moment problems is provided by another
result of Stochel (which actually holds for moment problems on $\mathbb{R}%
^{d}$, $d>1$).

\begin{theorem}
\label{thm14}(cf. \cite{Sto2}) \ $\beta^{(\infty)}$ has a representing measure
supported in a closed set $K\subseteq\mathbb{R}^{2}$ if and only if, for each
$n$, $\beta^{(2n)}$ has a representing measure supported in $K$.
\end{theorem}%

\noindent
In Section \ref{application} we will use Theorem \ref{thm14} to give a new
proof of Theorem \ref{thm13}. \ To do so, we require the following refinement
of Theorem \ref{firstmain}, which relates the existence of representing
measures to extensions of moment matrices.

\begin{theorem}
\label{hyper}Let $H:=\{(x,y)\in\mathbb{R}^{2}:Q(x,y)=0\}$ be an hyperbola.
\ For $\beta\equiv\beta^{(2n)}$, assume that $\mathcal{M}(n)\equiv
\mathcal{M}(n)(\beta)$ is positive, recursively generated, and satisfies
$Q(X,Y)=0$ in $\mathcal{C}_{\mathcal{M}(n)}$. \ Then $\operatorname*{rank}%
\mathcal{M}(n)\leq2n+1$, and the following statements are equivalent.
\end{theorem}

\begin{enumerate}
\item[(i)] $\beta$\textit{\ admits a representing measure (necessarily
supported in }$H$\textit{).}

\item[(ii)] $\beta$\textit{\ admits a representing measure (necessarily
supported in }$H$\textit{) with convergent moments of degree up to }$2n+2$\textit{.}

\item[(iii)] $\beta$\textit{\ admits a representing measure }$\mu
$\textit{\ (necessarily supported in }$H$\textit{) satisfying }%
$\operatorname*{card}\operatorname*{supp}\mu\leq1+\operatorname*{rank}%
M(n)$\textit{. \ If }$\operatorname*{rank}M(n)\leq2n$\textit{, then }$\mu
$\textit{\ can be taken so that }$\operatorname*{card}\operatorname*{supp}%
\mu=\operatorname*{rank}M(n)$\textit{.}

\item[(iv)] $M(n)$\textit{\ admits a positive, recursively generated moment
matrix extension }$M(n+1)$\textit{.}

\item[(v)] $M(n)$\textit{\ admits a positive, recursively generated extension
}$M(n+1)$\textit{, with }$\operatorname*{rank}M(n+1)\leq1+\operatorname*{rank}%
M(n)$\textit{, and }$M(n+1)$\textit{\ admits a flat extension }$M(n+2)$%
\textit{. \ If }$\operatorname*{rank}M(n)\leq2n$\textit{, then }%
$M(n)$\textit{\ admits a flat extension }$M(n+1)$\textit{.}

\item[(vi)] $\operatorname*{rank}\;M(n)\leq\operatorname*{card}\;V(M(n))$\textit{.}
\end{enumerate}

Condition (vi) in Theorem \ref{hyper} is the concrete condition which,
together with positivity and recursiveness, provides an effective test for the
existence of representing measures; we illustrate Theorem \ref{hyper} with an example.

\begin{example}
\label{ex16}We consider a case of $M(2)$ satisfying $YX=1$; let
\[
\mathcal{M}(2)(\beta):=\left(
\begin{array}
[c]{cccccc}%
1 & 0 & 0 & a & 1 & a\\
0 & a & 1 & 0 & 0 & 0\\
0 & 1 & a & 0 & 0 & 0\\
a & 0 & 0 & c & a & 1\\
1 & 0 & 0 & a & 1 & a\\
a & 0 & 0 & 1 & a & c
\end{array}
\right)  .
\]
By calculating nested determinants, we see at once that the $3$ by $3$ upper
left-hand corner is positive and invertible if and only if $a>1$, and that the
$4$ by $4$ upper left-hand corner is positive and invertible if and only if
$c>a^{2}$. \ We now let $a>1$ and $c:=a^{2}+r$, where $r>0$. \ Then
$\operatorname*{rank}M(2)=4$ if and only if
\begin{equation}
(a^{2}-1)(r^{2}-a^{2}+1)=0. \label{eq in r}%
\end{equation}
\ The positive root of (\ref{eq in r}) is $r_{1}:=\sqrt{a^{2}-1}$, so we set
$r:=r_{1}$ and observe that in $C_{\mathcal{M}(2)}$, $Y^{2}=2a1-X^{2}$.
\ Since $YX=1$, we next find the variety $V(\beta)$ by solving the pair of
equations
\begin{equation}
\left\{
\begin{array}
[c]{c}%
y^{2}=2a-x^{2}\\
y=\frac{1}{x}%
\end{array}
.\right.  \label{eq in x and y}%
\end{equation}
It is easy to see that (\ref{eq in x and y}) has exactly four roots,
$\{(x_{i},y_{i})\}_{i=1}^{4}$, where $x_{1}:=-(a-\sqrt{a^{2}-1})$,
$x_{2}:=a-\sqrt{a^{2}-1}$, $x_{3}:=-(a+\sqrt{a^{2}-1})$, $x_{4}:=a+\sqrt
{a^{2}-1}$, and $y_{i}:=\frac{1}{x_{i}}\;(i=1,2,3,4)$. \ Thus,
$\operatorname*{rank}M(2)=4=\operatorname*{card}V(M(2))$. \ According to
Theorem \ref{hyper}(vi)$\Rightarrow$(i)$,\beta^{(4)}$ admits a representing
measure $\mu$. \ Since $4=\operatorname*{rank}M(2)\leq\operatorname*{card}%
\operatorname*{supp}\mu\leq\operatorname*{card}V(\beta)=4$, it follows that
$\operatorname*{supp}\mu=V(\beta)$ and that $M(2)$ admits a flat extension
$M(3)$ (cf. (\ref{exp18})). \qed                
\end{example}

Theorem \ref{hyper} shows that minimal $Q$-representing measures for
$\beta^{(2n)}$ arise either from flat extensions of $\mathcal{M}(n)$ or of
$\mathcal{M}(n+1)$ (cf. \cite{FiaOT}). \ In the presence of a flat extension,
there is a simple procedure for computing the atoms and densities of a
corresponding minimal representing measure.

\begin{theorem}
\label{thm17}(cf. \cite[Theorem 2.21]{tcmp10}) \ If $\mathcal{M}%
(n)\equiv\mathcal{M}(n)(\beta)$ is positive semi-definite and admits a flat
extension $\mathcal{M}(n+1)$, then $\mathcal{V}:=\mathcal{V}(\mathcal{M}%
(n+1))$ satisfies $\operatorname*{card}\mathcal{V}=r\;(\equiv
\operatorname*{rank}\mathcal{M}(n))$, and $\mathcal{V}\equiv\{(x_{i}%
,y_{i})\}_{i=1}^{r}\subseteq\mathbb{R}^{2}$ forms the support of the unique
representing measure $\mu$ for $\mathcal{M}(n+1)$, i.e., $\mu$ is of the form
$\mu=\sum_{i=1}^{r}\rho_{i}\delta_{(x_{i},y_{i})}$ with $\rho_{i}>0\;(1\leq
i\leq r)$. \ If $\mathcal{B}\equiv\{Y^{i_{k}}X^{j_{k}}\}_{k=1}^{r}$ is a
maximal linearly independent subset of columns of $\mathcal{M}(n)$, let $V$ be
the $r\times r$ matrix whose entry in row $k$, column $\ell$ is $y_{\ell
}^{i_{k}}x_{\ell}^{j_{k}}\;(1\leq k,\ell\leq r)$. \ Then $V$ is invertible,
and $\rho\equiv(\rho_{1},...,\rho_{r})$ is uniquely determined by $V\rho
^{t}=(\beta_{i_{1},j_{1}},...,\beta_{i_{r},j_{r}})^{t}$. \ 
\end{theorem}

\begin{example}
(Example \ref{ex16} cont.) \ We now use Theorem \ref{thm17} to compute the
densities for the measure $\mu:=\sum_{i=1}^{4}\rho_{i}\delta_{(x_{i},y_{i})}$
of Example \ref{ex16}. \ Since $\{1,X,Y,X^{2}\}$ is a basis for $\mathcal{C}%
_{\mathcal{M}(2)}$, Theorem \ref{thm17} implies that
\[
V\equiv\left(
\begin{array}
[c]{cccc}%
1 & 1 & 1 & 1\\
x_{1} & x_{2} & x_{3} & x_{4}\\
y_{1} & y_{2} & y_{3} & y_{4}\\
x_{1}^{2} & x_{2}^{2} & x_{3}^{2} & x_{4}^{2}%
\end{array}
\right)  ,
\]
is invertible (indeed, $\det V=(x_{1}^{2}-x_{3}^{2})^{2}$). \ The vector
$\mathbf{v}:=(\beta_{i_{1},j_{1}},...,\beta_{i_{r},j_{r}})^{t}$ in Theorem
\ref{thm17} is $(1,0,0,a)^{t}$, so a calculation of $\rho=V^{-1}\mathbf{v}$
shows that $\rho_{1}=\rho_{2}=\rho_{3}=\rho_{4}=\frac{1}{4}$.
\textbf{\qed                }
\end{example}

To prove Theorems \ref{firstmain}, \ref{thm12} and \ref{hyper}, we will reduce
the analysis of truncated moment problems on conics to the study of truncated
moment problems on four special conics: $x^{2}+y^{2}=1$, $y=x^{2}$, $yx=1$,
and $yx=0$. \ This reduction was initially described (using complex moment
matrices) in \cite{tcmp6}, but we now require a more detailed analysis.

For $a_{1},a_{2},b_{1},b_{2},c_{1},c_{2}\in\mathbb{R}$, with $b_{1}c_{2}\neq
b_{2}c_{1}$, let $\Phi:\mathbb{R}^{2}\rightarrow\mathbb{R}^{2}$ denote the
degree-one bijection defined by $\Phi(x,y)\equiv(\varphi_{1},\varphi
_{2}):=(a_{1}+b_{1}x+c_{1}y,a_{2}+b_{2}x+c_{2}y)\;(x,y\in\mathbb{R)}$. \ Given
$\beta^{(2n)}$, define $\tilde{\beta}^{\left(  2n\right)  }$ by $\tilde{\beta
}_{ij}:=L_{\beta}(\varphi_{2}^{i}\varphi_{1}^{j})$ ($0\leq i+j\leq2n$), where
$L_{\beta}$ denotes the Riesz functional associated with $\beta$. It is
straightforward to verify that $L_{\tilde{\beta}}(p)=L_{\beta}\left(
p\circ\Phi\right)  $ for every $p\in\mathbb{R}_{n}[x,y]$. (Note that for
$p\left(  x,y\right)  \equiv\sum a_{ij}y^{i}x^{j}$, $\left(  p\circ
\Phi\right)  \left(  x,y\right)  =p\left(  \varphi_{1},\varphi_{2}\right)
\equiv\sum a_{ij}\varphi_{2}^{i}\varphi_{1}{}^{j}$.) \ Also, $L_{\beta
}(p)=L_{\tilde{\beta}}(p\circ\Psi)$, where $\Psi:=\Phi^{-1}$, given by
$\Psi(u,v):=\frac{1}{b_{1}c_{2}-b_{2}c_{1}}(c_{2}(u-a_{1})-c_{1}%
(v-a_{2}),-b_{2}(u-a_{1})+b_{1}(v-a_{2}))$.

Let $\mathcal{H}^{(n)}:=\{\hat{p}:p\in\mathbb{R}_{n}[x,y]\}\cong
\mathbb{R}^{m(n)}$ and consider the linear map $J\equiv J^{(n)}:\mathcal{H}%
^{(n)}\rightarrow\mathcal{H}^{(n)}$ defined by $J(\hat{p}):=\widehat
{p\circ\Phi}\;(p\in\mathbb{R}_{n}[x,y])$. \ The map $J$ is invertible, with
inverse $J^{-1}(\hat{q}):=\widehat{q\circ\Psi}\;(q\in\mathbb{R}_{n}[x,y])$.
\ Note that $\mathcal{H}^{(n)}$ admits a vector space direct sum decomposition
$\mathcal{H}^{(n)}=\mathcal{H}^{(n-1)}\bigoplus\mathcal{H}_{n} $, where
$\mathcal{H}_{n}$ is the subspace spanned by vectors $\widehat{y^{i}x^{j}}$
with $i+j=n$. \ Since $J^{(n)}$ and $J^{(n-1)}:\mathcal{H}^{(n-1)}%
\rightarrow\mathcal{H}^{(n-1)}$ are both invertible, relative to the above
decomposition $J^{(n)}$ can be represented as
\begin{equation}
J^{(n)}=\left(
\begin{array}
[c]{cc}%
J^{(n-1)} & K_{n}\\
0 & L_{n}%
\end{array}
\right)  , \label{eq110}%
\end{equation}
and thus $(J^{(n)})^{-1}$ admits a similar triangular representation.

\begin{proposition}
\label{prop17}\textup{(Invariance under degree-one transformations.)} \ Let
$\mathcal{M}(n)$ and $\mathcal{\tilde{M}}(n)$ be the moment matrices
associated with $\beta$ and $\tilde{\beta}$\textup{.}

\begin{enumerate}
\item[(i)] $\mathcal{\tilde{M}}(n)=J^{\ast}\mathcal{M}(n)J$.

\item[(ii)] $J$ is invertible.

\item[(iii)] $\mathcal{\tilde{M}}(n)\geq0\Leftrightarrow\mathcal{M}(n)\geq0$.

\item[(iv)] $\operatorname*{rank}\mathcal{\tilde{M}}(n)=\operatorname*{rank}%
\mathcal{M}(n)$.

\item[(v)] The formulas $\mu=\tilde{\mu}\circ\Phi$ and $\tilde{\mu}=\mu
\circ\Psi$ establish a one-to-one correspondence between the sets of
representing measures for $\beta$ and $\tilde{\beta}$, which preserves measure
class and cardinality of the support; moreover, $\Phi(\operatorname*{supp}%
\mu)=\operatorname*{supp}\tilde{\mu}$ and $\Psi(\operatorname*{supp}\tilde
{\mu})=\operatorname*{supp}\mu$.

\item[(vi)] For $p\in\mathbb{R}_{n}[x,y]$, $p\left(  \tilde{X},\tilde
{Y}\right)  =J^{\ast}\left(  \left(  p\circ\Phi\right)  \left(  X,Y\right)
\right)  $ and $p(X,Y)=$\newline $(J^{-1})^{\ast}((p\circ\Psi)(\tilde
{X},\tilde{Y}))$.

\item[(vii)] $\mathcal{\tilde{M}}(n)$ is recursively generated if and only if
$\mathcal{M}(n)$ is recursively generated.

\item[(viii)] $\Phi(\mathcal{V}(\beta))=\mathcal{V}(\tilde{\beta})$ and
$\mathcal{V}(\beta)=\Psi(\mathcal{V}(\tilde{\beta}))$.

\item[(ix)] $\mathcal{M}(n)$ is positive semi-definite and admits a positive,
recursively generated (resp. flat) extension $\mathcal{M}(n+1)$ if and only if
$\mathcal{\tilde{M}}(n)$ is positive semi-definite and admits a positive,
recursively generated (resp. flat) extension $\mathcal{\tilde{M}}(n+1)$.
\end{enumerate}
\end{proposition}

\begin{proof}
We omit the proofs of (i) - (viii), which are straightforward. \ For (ix),
suppose $\mathcal{M}(n)$ is positive and admits a positive, recursively
generated extension
\[
\mathcal{M}(n+1)(\check{\beta})\equiv\left(
\begin{array}
[c]{cc}%
\mathcal{M}(n) & B(n+1)\\
B(n+1)^{t} & C(n+1)
\end{array}
\right)  .
\]
\ Part (i) (using $n+1$ instead of $n$) and (\ref{eq110}) imply that
\begin{align*}
\mathcal{M}  &  \equiv\mathcal{\tilde{M}}(n+1)=(J^{(n+1)})^{\ast}%
\mathcal{M}(n+1)(\check{\beta})J^{(n+1)}\\
&  =\left(
\begin{array}
[c]{cc}%
(J^{(n)})^{\ast}\mathcal{M}(n)J^{(n)} & \ast\\
\ast & \ast
\end{array}
\right)  =\left(
\begin{array}
[c]{cc}%
\mathcal{\tilde{M}}(n) & \ast\\
\ast & \ast
\end{array}
\right)  .
\end{align*}
Thus, $\mathcal{M}$ is a positive moment matrix extension of $\mathcal{\tilde
{M}}(n)$, and (vii) (applied with $n+1$) implies that $\mathcal{M}$ is
recursively generated. \ In the case when $\mathcal{M}$ is a flat extension,
we have, from (i), (ii) and (iv), $\operatorname*{rank}\mathcal{\tilde{M}%
}(n)=\operatorname*{rank}\mathcal{M}(n)=\operatorname*{rank}\mathcal{M}%
(n+1)(\check{\beta})=\operatorname*{rank}\mathcal{M}$, so $\mathcal{M}$ is a
flat extension of $\mathcal{\tilde{M}}(n)$. \ The converse is proved in the
same way.
\end{proof}

It is well known that a general conic $Q(x,y)=0$ can be transformed into one
of the four special cases by means of rotation, scaling and translation; thus,
there is a degree-one map $\Phi$ (as above) such that $Q\circ\Phi$ is a
special conic. \ Note from Proposition \ref{prop17} that such properties of
$\mathcal{M}(n)$ as positivity, recursiveness, the variety condition, and the
existence of flat or positive extensions $\mathcal{M}(n+1) $ are invariant
under degree-one mappings, which also preserve the existence of representing
measures and the cardinality of the support. \ These observations show that
Theorems \ref{firstmain} and \ref{hyper} are valid for arbitrary hyperbolas if
and only if they are valid for $yx=1$ and $yx=0$. \ We prove Theorem
\ref{firstmain} for $yx=1$ in Section \ref{3tm} (Theorem \ref{thmflat}) and
for $yx=0$ in Section \ref{degenerate} (Theorem \ref{secondflat}). \ We prove
Theorem \ref{hyper} (for both hyperbolas) in Section \ref{proof} (Theorem
\ref{thmhyperbolic2}). \ Section \ref{examples} contains examples illustrating
our results. \ Section \ref{application} contains a new proof of Theorem
\ref{thm13}, based on moment matrix techniques, including Theorem \ref{hyper}.

The remainder of this section is devoted to notation and basic results
concerning real moment matrices. \ Given a collection $\beta^{(2n)}:\beta
_{00},\beta_{01},\beta_{10},...,\beta_{0,2n},...,\beta_{2n,0}$, we can
describe $\mathcal{M}(n)(\beta)$ by means of a block matrix decomposition
$\mathcal{M}(n)(\beta):=(\mathcal{M}_{ij}(\beta))_{i,j=0}^{n}$, where
\[
\mathcal{M}_{ij}(\beta):=\left(
\begin{array}
[c]{cccc}%
\beta_{0,i+j} & \beta_{1,i+j-1} & \cdots & \beta_{j,i}\\
\beta_{1,i+j-1} & \beta_{2,i+j-2} & \cdots & \beta_{j+1,i-1}\\
\vdots & \vdots & \ddots & \vdots\\
\beta_{i,j} & \beta_{i+1,j-1} & \cdots & \beta_{i+j,0}%
\end{array}
\right)  .
\]
Recall that $\mathcal{M}(n)$ has size $m(n)\equiv\frac{(n+1)(n+2)}{2}$.
$\ $For any matrix $M$ of this size, $\left[  M\right]  _{k}$ denotes the
compression of $M$ to the first $k$ rows and columns; similarly, for a vector
$\mathbf{v}$, $\left[  \mathbf{v}\right]  _{k}$ denotes the compression of
$\mathbf{v}$ to the first $k$ entries. \ We also consider compressions of $M$
and $\mathbf{v}$ to a set $E$ of rows and columns, and denote such
compressions by $\left[  M\right]  _{E}$ and $\left[  \mathbf{v}\right]  _{E}%
$. \ For $i+j,k+\ell\leq n$, $\left\langle Y^{i}X^{j},Y^{k}X^{l}\right\rangle
_{M}$ (or simply $\left\langle Y^{i}X^{j},Y^{k}X^{l}\right\rangle $) denotes
the entry of $M$ in row $Y^{k}X^{l}$ and column $Y^{i}X^{j}$.

A theorem of Smul'jan \cite{Smu} shows that a block matrix
\begin{equation}
M=%
\begin{pmatrix}
A & B\\
B^{\ast} & C
\end{pmatrix}
\label{locaeqFlat.3}%
\end{equation}
is positive if and only if (i) $A\geq0$, (ii) there exists a matrix $W$ such
that $B=AW$, and (iii) $C\geq W^{\ast}AW$ (since $A=A^{\ast}$, $W^{\ast}AW$ is
independent of $W$ satisfying $B=AW$). \ Note also that if $M\geq0$, then
$\operatorname*{rank}M=\operatorname*{rank}A$ if and only if $C=W^{\ast}AW$;
conversely, if $A\geq0$ and there exists $W$ such that $B=AW$ and $C=W^{\ast
}AW$, then $M\geq0$ and $\operatorname*{rank}M=\operatorname*{rank}A$. \ In
the sequel, for $A\geq0$, we refer to $M$ as an \emph{extension} of $A$, and
as a \emph{flat extension} if $\operatorname*{rank}M=\operatorname*{rank}A$.
\ Thus, a flat extension of a positive matrix $A$ is completely determined by
a choice of block $B$ satisfying $B=AW$ and $C=W^{\ast}AW$ for some matrix
$W$; we denote such an extension by $\left[  A;B\right]  $. \ It follows from
the Extension Principle \cite{FiCM} that if $M\geq0$, then column dependence
relations in $A$ extend to $\left(
\begin{array}
[c]{c}%
A\\
B^{\ast}%
\end{array}
\right)  $; further, if $M$ is a flat extension of $A\;(\geq0)$, then column
dependence relations in $\left(  A\;\;B\right)  $ extend to $M$.

For an $\left(  n+1\right)  \times\left(  n+2\right)  $ moment matrix block
$B_{n,n+1}$, representing ``new moments'' of degree $2n+1$ for a prospective
representing measure of $\beta^{\left(  2n\right)  }$, let
\[
B(n+1):=%
\begin{pmatrix}
B_{0,n+1}\\
\vdots\\
B_{n-1,n+1}\\
B_{n,n+1}%
\end{pmatrix}
.
\]
By Smul'jan's theorem, $\mathcal{M}(n)\geq0$ admits a (necessarily positive)
flat extension
\[
\left[  \mathcal{M}(n);B\right]  =%
\begin{pmatrix}
\mathcal{M}(n) & B\\
B^{\ast} & C
\end{pmatrix}
\]
in the form of a moment matrix $\mathcal{M}(n+1)$ if and only if
\begin{equation}%
\begin{array}
[t]{l}%
B=B(n+1)\text{ and }B=\mathcal{M}(n)W\text{ for some }W\text{ }\\
\text{(i.e., Ran }B\subseteq\;\text{Ran }\mathcal{M}(n)\text{\cite{Dou});
and}\\
C:=W^{\ast}\mathcal{M}(n)W\text{ is Hankel}\\
\text{(i.e., }C\text{ has the form of a moment matrix block }B_{n+1,n+1}%
\text{).}%
\end{array}
\label{locaeqFlat.4}%
\end{equation}

\textit{Acknowledgement}. \ The examples in this paper were obtained using
calculations with the software tool \textit{Mathematica \cite{Wol}}.

\section{\label{3tm}The Truncated Moment Problem on Nondegenerate Hyperbolas}

In this section we prove Theorem \ref{firstmain} for nondegenerate hyperbolas.
\ In view of Proposition \ref{prop17} and the remarks following it, it
suffices to consider the case $yx=1$. \ The necessity of the conditions in
Theorem \ref{firstmain} is clear from Section \ref{Int}, and sufficiency
follows immediately from Theorem \ref{thmflat} below.

\begin{theorem}
\label{thmflat}Let $\beta\equiv\beta^{(2n)}:\beta_{00},\beta_{01},\beta
_{10},...,\beta_{0,2n},...,\beta_{2n,0}$ be a family of real numbers,
$\beta_{00}>0$, and let $\mathcal{M}(n)$ be the associated moment matrix.
$\ $Assume that $\mathcal{M}(n)$ is positive, recursively generated, and
satisfies $YX=1$ and $\operatorname*{rank}\;\mathcal{M}(n)\leq
\operatorname*{card}\;\mathcal{V}(\beta)$. \ Then $\operatorname*{rank}%
\mathcal{M}(n)\leq2n+1$. \ If $\operatorname*{rank}\mathcal{M}(n)\leq2n$, then
$\mathcal{M}(n)$ admits a flat extension $\mathcal{M}(n+1)$ (so $\beta$ admits
a $\operatorname*{rank}\mathcal{M}(n)$-atomic representing measure supported
in $yx=1$). \ If $\operatorname*{rank}\mathcal{M}(n)=2n+1$, then
$\mathcal{M}(n)$ admits an extension to a positive, recursively generated
extension $\mathcal{M}(n+1)$, satisfying $2n+1\leq\operatorname*{rank}%
\mathcal{M}(n+1)\leq2n+2$, and $\mathcal{M}(n+1)$ admits a flat extension
$\mathcal{M}(n+2)$ (so $\beta$ admits a representing measure $\mu$ supported
in $yx=1$, with $2n+1\leq\operatorname*{card}\operatorname*{supp}\mu\leq2n+2$.
\end{theorem}

We require several preliminary results for the proof of Theorem \ref{thmflat}.
\ By \cite[Theorem 2.1]{tcmp2}, we know that $\mathcal{M}(n)$ (positive
semi-definite and recursively generated) admits flat extensions when $\left\{
\mathit{1},X,Y\right\}  $ is linearly dependent in $\mathcal{C}_{\mathcal{M}%
(n)}$. \ Thus, hereafter we will assume that $\left\{  \mathit{1},X,Y\right\}
$ is linearly independent. \ We begin with an elementary lemma that exploits
the fact that $\mathcal{M}(n)$ is recursively generated. \ For $1\leq k\leq n$
let
\[
\mathcal{S}_{n}(k):=\{1,X,Y,X^{2},Y^{2},...,X^{k},Y^{k}\}\subseteq
\mathcal{C}_{\mathcal{M}(n)}.
\]

\begin{lemma}
\label{lem2.3}For $n\geq2$, let $\mathcal{M}(n)$ be positive and recursively
generated, and assume that $YX=1$. \ Then each column of $\mathcal{M}(n)$ is
equal to a column in $\mathcal{S}_{n}(n)$; in particular,
$\operatorname*{rank}\mathcal{M}(n)\leq2n+1$.
\end{lemma}

\begin{proof}
The proof is by induction on $n\geq2$. \ For $n=2$ the statement is clearly
true, so assume it holds for $n=k\;(\geq2)$. \ Suppose $\mathcal{M}\left(
k+1\right)  $ is positive and recursively generated, with $YX=1$ in
$\mathcal{C}_{\mathcal{M}\left(  k+1\right)  }$. \ Let $i,j\geq0,\;i+j\leq k$.
\ By the induction hypothesis, each column of the form $[Y^{i}X^{j}]_{m(k)}$
is in $\mathcal{S}_{k}(k)$, and since $\mathcal{M}(k+1)\geq0$, the Extension
Principle \cite[Proposition 2.4]{FiCM} shows that $Y^{i}X^{j}\in
\mathcal{S}_{k}(k+1)\;(\subseteq\mathcal{S}_{k+1}(k+1))$. \ Since
$X^{k+1},Y^{k+1}\in\mathcal{S}_{k+1}(k+1)$, it now suffices to consider a
column in $\mathcal{M}(k+1)$ of the form $Y^{k+1-j}X^{j}$, with $1\leq j\leq
k$. \ Let $q(x,y):=yx-1$ and let $p_{ij}(x,y):=y^{i}x^{j}$, so that
$Y^{k+1-j}X^{j}=p_{k+1-j,j}(X,Y)$. \ Also, let $r_{ij}(x,y):=y^{k+1-j}%
x^{j}-y^{k-j}x^{j-1}$. \ Now $r_{ij}(x,y)=y^{k-j}x^{j-1}(yx-1)=p_{k-j,j-1}%
(x,y)q(x,y)$; since $\mathcal{M}\left(  k+1\right)  $ is recursively generated
and $q(X,Y)=0$, it follows that $r_{ij}(X,Y)=0$, that is, $Y^{k+1-j}%
X^{j}=Y^{k-j}X^{j-1}$ in $\mathcal{C}_{\mathcal{M}(k+1)}$. \ By induction,
$[Y^{k-j}X^{j-1}]_{m(k)}\in\mathcal{S}_{k}(k)$, and since $\mathcal{M}\left(
k+1\right)  \geq0$, it follows as above that $Y^{k-j}X^{j-1}\in\mathcal{S}%
_{k+1}(k)$. \ Thus $Y^{k+1-j}X^{j}\;(=Y^{k-j}X^{j-1})\in\mathcal{S}%
_{k+1}(k)\subseteq\mathcal{S}_{k+1}(k+1)$, as desired.
\end{proof}

We next present two auxiliary results that will be used frequently in the
sequel. \ Recall that for $i+j,k+\ell\leq n$, $\left\langle Y^{i}X^{j}%
,Y^{k}X^{\ell}\right\rangle $ denotes the entry of $\mathcal{M}(n) $ in row
$Y^{k}X^{\ell}$, column $Y^{i}X^{j}$, namely $\beta_{i+k,j+\ell}$. \ We extend
this inner product notation from monomials to polynomials as follows. \ For
$p\equiv\sum_{0\leq i+j\leq n}a_{ij}y^{i}x^{j} $ and $q\equiv\sum_{0\leq
k+\ell\leq n}b_{k\ell}y^{k}x^{\ell}$, we define\ $\left\langle
p(X,Y),q(X,Y)\right\rangle :=\sum_{0\leq i+j,k+\ell\leq n}a_{ij}b_{k\ell}%
\beta_{i+k,j+\ell}$. \ Further, if $\deg p+\deg p^{\prime},\deg q+\deg
q^{\prime}\leq n$, by\ $\left\langle p(X,Y)p^{\prime}(X,Y),q(X,Y)q^{\prime
}(X,Y)\right\rangle $ we mean $\left\langle (pp^{\prime})(X,Y),(qq^{\prime
})(X,Y)\right\rangle $. \ The following result follows directly from the
preceding definitions.

\begin{lemma}
\label{newlem}(i) \ For $p,q\in\mathbb{R}_{n}[x,y]$,
\[
\left\langle p(X,Y),q(X,Y)\right\rangle =\left\langle
q(X,Y),p(X,Y)\right\rangle .
\]
\newline (ii) \ For $p,q\in\mathbb{R}_{n}[x,y]$, $i,j\geq0$, $i+j\leq n$, and
$\deg p,\deg q\leq n-(i+j)$,
\[
\left\langle p(X,Y)Y^{j}X^{i},q(X,Y)\right\rangle =\left\langle
p(X,Y),q(X,Y)Y^{j}X^{i}\right\rangle .
\]
\newline (iii) \ If $p,q,r\in\mathbb{R}_{n}[x,y]$ with $p(X,Y)=q(X,Y)$ in
$\mathcal{C}_{\mathcal{M}(n)}$, then $\left\langle r(X,Y),p(X,Y)\right\rangle
$\newline $=\left\langle r(X,Y),q(X,Y)\right\rangle $.
\end{lemma}

\begin{lemma}
\label{basic}Let $\mathcal{M}(n)$ be positive, recursively generated, with
$YX=1$, and assume $p,q\in\mathbb{R}_{n-1}[x,y]$. \ Then
\begin{align}
\left\langle p(X,Y),q(X,Y)\right\rangle  &  =\left\langle
Yp(X,Y),Xq(X,Y)\right\rangle \label{eq101}\\
&  =\left\langle Xp(X,Y),Yq(X,Y)\right\rangle . \label{eq102}%
\end{align}
\end{lemma}

\begin{proof}
The definition of $\left\langle p(X,Y),q(X,Y)\right\rangle $ implies that,
without loss of generality, we can assume that $p(X,Y)=Y^{i}X^{j}$ and
$q(X,Y)=Y^{k}X^{\ell}$. \ Assume first that $k\geq1$. \ We have
\begin{align*}
\left\langle Y^{i}X^{j},Y^{k}X^{\ell}\right\rangle  &  =\left\langle
Y^{i+1}X^{j},Y^{k-1}X^{\ell}\right\rangle \;\;\text{(by Lemma \ref{newlem}%
(ii))}\\
&  =\left\langle Y^{i+1}X^{j},Y^{k-1}X^{\ell}YX\right\rangle \;\;\text{(by
Lemma \ref{newlem}(iii), using }YX=1\text{)}\\
&  =\left\langle Y^{i+1}X^{j},Y^{k}X^{\ell+1}\right\rangle .
\end{align*}
If $k=0$ and $j\geq1$, we have
\begin{align*}
\left\langle Y^{i}X^{j},X^{\ell}\right\rangle  &  =\left\langle Y^{i}%
X^{j-1},X^{\ell+1}\right\rangle \;\;\text{(by Lemma \ref{newlem}(ii))}\\
&  =\left\langle Y^{i}X^{j-1}YX,X^{\ell+1}\right\rangle \;\;\text{(by Lemma
\ref{newlem}(iii), using }YX=1\text{)}\\
&  =\left\langle Y^{i+1}X^{j},X^{\ell+1}\right\rangle .
\end{align*}
If $k=j=0$, we have $p(X,Y)=Y^{i}$, $q(X,Y)=X^{\ell}$, so we need to prove
that $\left\langle Y^{i},X^{\ell}\right\rangle =\left\langle Y^{i+1}%
,X^{\ell+1}\right\rangle $. \ If $i\geq1$, we have
\begin{align*}
\left\langle Y^{i},X^{\ell}\right\rangle  &  =\left\langle Y^{i-1},YX^{\ell
}\right\rangle \;\;\text{(by Lemma \ref{newlem}(ii))}\\
&  =\left\langle Y^{i}X,YX^{\ell}\right\rangle \;\;\text{(by Lemma
\ref{newlem}(i),(iii))}\\
&  =\beta_{i+1,\ell+1}=\left\langle Y^{i+1},X^{\ell+1}\right\rangle \text{.}%
\end{align*}
If $i=0$, then
\begin{align*}
\left\langle Y^{i},X^{\ell}\right\rangle  &  =\left\langle 1,X^{\ell
}\right\rangle =\left\langle YX,X^{\ell}\right\rangle \;\;\text{(by Lemma
\ref{newlem}(iii))}\\
&  \left\langle Y,X^{\ell+1}\right\rangle \;\;\text{(by Lemma \ref{newlem}%
(ii)).}%
\end{align*}
We have now completed the proof of (\ref{eq101}); the validity of
(\ref{eq102}) is a straightforward consequence of (\ref{eq101}) and Lemma
\ref{newlem}(i).
\end{proof}

We next divide the proof of Theorem \ref{thmflat} into four cases, based on
possible dependence relations among the elements of $\mathcal{S}_{n}(n)$.
\ Section \ref{examples} contains examples illustrating these cases. \ In
proving each case, we ultimately obtain some flat moment matrix extension
$\mathcal{M}$; the existence of a corresponding $\operatorname*{rank}%
\mathcal{M}$-atomic representing measure $\mu$ supported in $yx=1$ then
follows immediately from (\ref{exp18}) and (\ref{15}); for this reason, and to
simplify the statement of each case, we address only the matrix extension, not
the representing measure. \ In the sequel, unless otherwise noted, we are
always assuming that $\mathcal{M}(n)$ is positive, recursively generated,
$\left\{  1,X,Y\right\}  $ is linearly independent, $YX=1$, and
$\operatorname*{rank}\mathcal{M}(n)\leq\operatorname*{card}\mathcal{V}(\beta
)$. \ 

\begin{proposition}
\label{propcaseI}(Case I: For some $k$, $2\leq k\leq n$, $\mathcal{S}%
_{n}(k-1)$ is linearly independent and $X^{k}\in\;$lin.span $\mathcal{S}%
_{n}(k-1)$) \ Assume that $\mathcal{M}(n)(\beta)$ is positive, recursively
generated, satisfies $YX=1$, and $\operatorname*{rank}\mathcal{M}%
(n)\leq\operatorname*{card}\mathcal{V}(\beta)$. \ In $\mathcal{S}_{n}(n)$,
assume that the first dependence relation occurs at $X^{k}$, with $2\leq k\leq
n$. \ Then $\mathcal{M}(n)$ is flat and, a fortiori, it admits a unique flat
extension $\mathcal{M}(n+1)$.
\end{proposition}

\begin{proof}
Write $X^{k}=p_{k-1}(X)+q_{k-1}(Y),$where $\deg p_{k-1},\deg q_{k-1}\leq k-1$.
\ It follows that $\mathcal{V}(\beta)\subseteq(yx=1)%
{\textstyle\bigcap}
(p_{k-1}(x)+q_{k-1}(y)=x^{k})$ $\subseteq(yx=1)%
{\textstyle\bigcap}
(p_{k-1}(x)+q_{k-1}(\frac{1}{x})=x^{k})$. $\ $Since $p_{k-1}(x)+q_{k-1}%
(\frac{1}{x})=x^{k}$ leads to a polynomial equation in $x$ of degree at most
$2k-1$, it follows that $\operatorname*{card}\;\mathcal{V}(\beta)\leq2k-1$, so
$\operatorname*{rank}\mathcal{M}(n)\leq2k-1$. \ Then $\mathcal{S}%
_{n}(k-1)\equiv\{1,X,Y,X^{2},Y^{2},...,X^{k-1},Y^{k-1}\}$ is a basis for
$\mathcal{C}_{\mathcal{M}(n)}$, whence $\mathcal{M}(n)$ is flat.
\end{proof}

\begin{proposition}
\label{propcaseIII}(Case II: For some $k$, $2\leq k<n$, $\mathcal{S}%
_{n}(k-1)\bigcup\{X^{k}\}$ is linearly independent, and $Y^{k}\in\;$lin.span
$(\mathcal{S}_{n}(k-1)\bigcup\{X^{k}\}$) \ Assume that $\mathcal{M}(n)(\beta)$
is positive, recursively generated, and satisfies $YX=1$. \ In $\mathcal{S}%
_{n}(n)$, assume that the first dependence relation occurs at $Y^{k}$, with
$1\leq k<n$. \ Then $\mathcal{M}(n)$ is flat, and thus admits a unique flat
extension $\mathcal{M}(n+1)$.
\end{proposition}

\begin{proof}
Write $Y^{k}=p_{k}(X)+q_{k-1}(Y)$, where $\deg p_{k}\leq k$ and $\deg
q_{k-1}\leq k-1$. \ Since $Y^{k}$ corresponds to a monomial of degree at most
$n-1$, and since $YX=1$ and $\mathcal{M}(n)$ is recursively generated, we must
have
\begin{equation}
Y^{k-1}=XY^{k}=Xp_{k}(X)+Xq_{k-1}(Y).\ \label{eq203}%
\end{equation}
Since $\mathcal{M}(n)$ is recursively generated and $YX=1$, $Xq_{k-1}(Y)$ is
clearly a linear combination of columns corresponding to monomials of degree
at most $k-2$. \ Let $a_{k}$ be the coefficient of $X^{k}$ in $p_{k}$. \ If
$a_{k}=0$, it follows from (\ref{eq203}) that $\mathcal{S}_{n}(k-1)\bigcup
\{X^{k}\}$ is linearly dependent, a contradiction. \ Thus, we must have
$a_{k}\neq0$, whence (\ref{eq203}) implies that $X^{k+1}$ is a linear
combination of previous columns. \ Moreover, $Y^{k+1}=Yp_{k}(X)+Yq_{k-1}(Y)$,
and $Yp_{k}(X)$ has degree $k-1$ in $X$, so $\mathcal{M}(k+1)$ is flat. \ It
now follows from the Extension Principle \cite{FiCM} and recursiveness that
$\mathcal{M}(n)$ is flat, i.e., $\mathcal{M}(n)$ is a flat extension of
$\mathcal{M}(k)$.
\end{proof}

\begin{proposition}
\label{propcaseIV}(Case III: The first dependence relation occurs at $Y^{n}%
$)\ \ Assume that $\mathcal{M}(n)(\beta)$ is positive, recursively generated,
satisfies $YX=1$, and $\operatorname*{rank}\mathcal{M}(n)\leq
\operatorname*{card}\mathcal{V}(\beta)$. \ In $\mathcal{S}_{n}(n)$, assume
that $Y^{n}$ is the location of the first dependence relation. \ Then
$\mathcal{M}(n)$ admits a flat extension $\mathcal{M}(n+1)$.
\end{proposition}

The proof of Proposition \ref{propcaseIV} will require several preliminary
results (Lemmas \ref{lem1}-\ref{yn+1xk} below). \ Under the hypotheses of
Proposition \ref{propcaseIV}, write
\begin{equation}
Y^{n}=a_{n}X^{n}+p_{n-1}(X)+q_{n-1}(Y), \label{yneq2}%
\end{equation}
with $\deg p_{n-1},\deg q_{n-1}\leq n-1$. \ We claim that $a_{n}\neq0$.
\ Assume instead that $a_{n}=0$, i.e., $Y^{n}=p_{n-1}(X)+q_{n-1}(Y)$. \ Then
$\mathcal{V}(\beta)\subseteq(yx=1)%
{\textstyle\bigcap}
(p_{n-1}(x)+q_{n-1}(y)=y^{n})$ $\subseteq(yx=1)%
{\textstyle\bigcap}
(p_{n-1}(\frac{1}{y})+q_{n-1}(y)=y^{n})$. $\ $Since $p_{n-1}(\frac{1}%
{y})+q_{n-1}(y)=y^{n}$ leads to a polynomial equation in $y$ of degree at most
$2n-1$, it follows that $\operatorname*{card}\;\mathcal{V}(\beta)\leq2n-1$, so
$\operatorname*{rank}\mathcal{M}(n)\leq2n-1$. \ Then $\mathcal{S}%
_{n}(n-1)\equiv\{1,X,Y,X^{2},Y^{2},...,X^{n-1},Y^{n-1}\}$ is a basis for
$\mathcal{C}_{\mathcal{M}(n)}$, whence $X^{n}$ is a linear combination of the
columns in $\mathcal{S}_{n}(n-1)$, a contradiction. \ Thus, $a_{n}\neq0$, so
in particular
\begin{equation}
X^{n}=\frac{1}{a_{n}}[Y^{n}-p_{n-1}(X)-q_{n-1}(Y)]. \label{xn}%
\end{equation}
To build a flat extension $\mathcal{M}(n+1)\equiv\left(
\begin{array}
[c]{cc}%
\mathcal{M}(n) & B(n+1)\\
B(n+1)^{\ast} & C(n+1)
\end{array}
\right)  $, we define the middle $n$ columns of a prospective block $B\equiv
B(n+1)$ by exploiting recursiveness and the relation $YX=1$, as follows:
\begin{equation}
YX^{n}:=X^{n-1};\;Y^{2}X^{n-1}:=YX^{n-2};\;...,\;Y^{n}X:=Y^{n-1}%
.\ \label{intermediate}%
\end{equation}
Also, motivated by (\ref{xn}) and, respectively, (\ref{yneq2}), we let
\begin{equation}
X^{n+1}:=\frac{1}{a_{n}}[Y^{n-1}-Xp_{n-1}(X)-Xq_{n-1}(Y)], \label{xn+1}%
\end{equation}
and
\begin{equation}
Y^{n+1}:=a_{n}X^{n-1}+Yp_{n-1}(X)+Yq_{n-1}(Y). \label{yn+1}%
\end{equation}
\ (The expressions $Y^{n-1}-Xp_{n-1}(X)-Xq_{n-1}(Y)$ and $a_{n}X^{n-1}%
+Yp_{n-1}(X)+Yq_{n-1}(Y)$ are shorthand notation for $(y^{n-1}-xp_{n-1}%
(x)-xq_{n-1}(y))(X,Y)$ and $(a_{n}x^{n-1}+yp_{n-1}(x)+yq_{n-1}(y))(X,Y)$ in
$\mathcal{C}_{\mathcal{M}(n)} $, respectively. \ Observe that these defining
relations are all required if one is to obtain a positive recursively
generated moment matrix extension $\mathcal{M}(n+1)$.) \ Since the columns
defined by (\ref{intermediate}) - (\ref{yn+1}) belong to $\mathcal{C}%
_{\mathcal{M}(n)}$, we have $B=\mathcal{M}(n)W$ for some matrix $W$. \ Thus, a
flat extension $M:=[\mathcal{M}(n);B]$ is uniquely determined by defining the
$C$-block as $C:=W^{\ast}\mathcal{M}(n)W$ (cf. Section 1). \ To complete the
proof that $M $ is a moment matrix $\mathcal{M}(n+1)$, it suffices to show
that block $B$ is of the form $(B_{i,n+1})_{i=0}^{n}$ and that block $C$ is of
the form $B_{n+1,n+1}$. \ To this end, we require some additional notation and
several preliminary results.

We next extend the notation $\left\langle p(X,Y),q(X,Y)\right\rangle $ to the
case when $\deg p=n+1$, $\deg q\leq n$. \ Indeed, using the definitions of the
columns of $B$, for $i,j\geq0$, $i+j=n+1$, there exists $p_{ij}\in
\mathbb{R}_{n}[x,y]$ with $Y^{i}X^{j}=p_{ij}(X,Y)$, and we define
\[
\left\langle Y^{i}X^{j},q(X,Y)\right\rangle :=\left\langle p_{ij}%
(X,Y),q(X,Y)\right\rangle .\
\]
Now if $p(x,y)\equiv\sum_{0\leq k+\ell\leq n+1}a_{k\ell}x^{\ell}y^{k}$, we
define
\[
\left\langle p(X,Y),q(X,Y)\right\rangle :=\sum_{0\leq k+\ell\leq n+1}a_{k\ell
}\left\langle Y^{k}X^{\ell},q(X,Y)\right\rangle .
\]

It is easy to check that Lemma \ref{newlem}(iii) holds with $\deg r=n+1$.

\begin{lemma}
\label{lem1}Under the hypotheses of Proposition \ref{propcaseIV}, assume
$i,j\geq0$, with $i+j=n+1$, and $r,s\geq1$, with $r+s\leq n$. \ Then
\begin{equation}
\left\langle Y^{i}X^{j},Y^{r}X^{s}\right\rangle =\left\langle Y^{i}%
X^{j},Y^{r-1}X^{s-1}\right\rangle . \label{eqlem1}%
\end{equation}
\end{lemma}

\begin{proof}
Fix $i$ and $j$ with $i+j=n+1$. \ We know from (\ref{intermediate}) -
(\ref{yn+1}) that there exists a polynomial $p\in\mathbb{R}_{n}[x,y]$ such
that $Y^{i}X^{j}=p(X,Y)\equiv\sum_{k+\ell\leq n}a_{k,\ell}Y^{k}X^{\ell}$.
\ Then
\begin{align*}
\left\langle Y^{i}X^{j},Y^{r}X^{s}\right\rangle  &  =\sum_{k+\ell\leq
n}a_{k,\ell}\left\langle Y^{k}X^{\ell},Y^{r}X^{s}\right\rangle \\
&  =\sum_{k+\ell\leq n}a_{k,\ell}\left\langle Y^{r}X^{s},Y^{k}X^{\ell
}\right\rangle \text{ \ (because }\mathcal{M}(n)\text{ is self-adjoint)}\\
&  =\sum_{k+\ell\leq n}a_{k,\ell}\left\langle Y^{r-1}X^{s-1},Y^{k}X^{\ell
}\right\rangle \;\;\text{(using }YX=1\text{ and recursiveness)}\\
&  =\sum_{k+\ell\leq n}a_{k,\ell}\left\langle Y^{k}X^{\ell},Y^{r-1}%
X^{s-1}\right\rangle \text{ \ }\\
&  \text{(using again the self-adjointness of }\mathcal{M}(n)\text{)}\\
&  =\left\langle Y^{i}X^{j},Y^{r-1}X^{s-1}\right\rangle \text{,}%
\end{align*}
as desired.
\end{proof}

The next result provides a reduction for the proof that $B(n+1)$ has the
Hankel property.

\begin{lemma}
\label{cor2}Under the hypotheses of Proposition \ref{propcaseIV}, assume
$i+j=n+1$, with $j\geq1$, $i\geq0$, and assume that the Hankel property
\begin{equation}
\left\langle Y^{i}X^{j},Y^{r}X^{s}\right\rangle =\left\langle Y^{i+1}%
X^{j-1},Y^{r-1}X^{s+1}\right\rangle \label{eqcor2}%
\end{equation}
holds with $1\leq r\leq n$ and $s=0$. $\ $Then (\ref{eqcor2}) holds for all
$r$ and $s$ such that $1\leq r+s\leq n$, $r\geq1$, $s\geq0$.
\end{lemma}

\begin{proof}
Fix $i$ and $j$ with $i+j=n+1$. \ We use induction on $t:=r+s$, where $1\leq
r+s\leq n$, $r\geq1$, $s\geq0$. \ For $t=1$ the result follows from the
hypothesis, since $r=1,$ $s=0$. $\ $Assume now that $t=2$. \ By hypothesis, we
may assume $r=s=1$, so we consider the equation
\begin{equation}
\left\langle Y^{i}X^{j},YX\right\rangle =\left\langle Y^{i+1}X^{j-1}%
,X^{2}\right\rangle , \label{eq211}%
\end{equation}
with $j\geq1$, $i\geq0$, $i+j=n+1$. \ Since $YX=1$, the left-hand side of
(\ref{eq211}) equals $\left\langle Y^{i}X^{j},1\right\rangle $ by Lemma
\ref{lem1}. \ For $j\geq2$ and $i\geq1$, the right-hand side of (\ref{eq211})
equals $\left\langle Y^{i}X^{j-2},X^{2}\right\rangle $ (by (\ref{intermediate}%
)), which in turn equals $\left\langle Y^{i}XX^{j-2},X\right\rangle
=\left\langle Y^{i-1}X^{j-2},X\right\rangle =\left\langle Y^{i-1}%
X^{j-1},1\right\rangle =\left\langle Y^{i}X^{j},1\right\rangle $ (by
(\ref{intermediate}) for the last step). \ When $j\geq2$ and $i=0$ (which then
implies $j=n+1)$, we have
\begin{align*}
\left\langle X^{n+1},YX\right\rangle  &  =\left\langle \frac{1}{a_{n}}%
[Y^{n-1}-Xp_{n-1}(X)-Xq_{n-1}(Y)],YX\right\rangle \text{ (by (\ref{xn+1})}\\
&  =\frac{1}{a_{n}}[\left\langle Y^{n},X\right\rangle -\left\langle
p_{n-1}(X),XYX\right\rangle -\left\langle q_{n-1}(Y),XYX\right\rangle ]\text{
}\\
&  \text{(by Lemma \ref{newlem}(ii) for the first term and Lemma \ref{basic}
for the }\\
&  \text{last two terms)}\\
&  =\left\langle \frac{1}{a_{n}}[Y^{n}-p_{n-1}(X)-q_{n-1}(Y)],X\right\rangle
\text{ (by Lemma \ref{newlem}(iii))}\\
&  =\left\langle X^{n},X\right\rangle \text{ (by (\ref{xn}))}\\
&  =\left\langle X^{n-1},X^{2}\right\rangle \\
&  =\left\langle YX^{n},X^{2}\right\rangle \text{ (by \ref{intermediate})).}%
\end{align*}
When $j=1$ (so that $i=n$), the right-hand side of (\ref{eq211}) is
\begin{align*}
\left\langle Y^{n+1},X^{2}\right\rangle  &  =\left\langle a_{n}X^{n-1}%
+Yp_{n-1}(X)+Yq_{n-1}(Y),X^{2}\right\rangle \;\;\text{(by(\ref{yn+1}))}\\
&  =\left\langle a_{n}X^{n}+p_{n-1}(X)+q_{n-1}(Y),X\right\rangle \;\;\\
&  \text{(using Lemma \ref{newlem}(ii) and Lemma \ref{basic}, as above)}\\
&  =\left\langle Y^{n},X\right\rangle \;\;\text{(by (\ref{yneq2}))}\\
&  =\left\langle Y^{n-1},1\right\rangle \;\;\text{(again using Lemma
\ref{basic}).}%
\end{align*}
On the other hand, the left-hand side of (\ref{eq211}) is
\begin{align*}
\left\langle Y^{n}X,YX\right\rangle  &  =\left\langle Y^{n-1},YX\right\rangle
\;\text{(by(\ref{intermediate}))}\\
&  =\left\langle Y^{n-1},1\right\rangle \text{ \ (by Lemma \ref{newlem}(iii)}.
\end{align*}
This completes the case when $t=2$. $\ $

Assume now that (\ref{eqcor2}) is true for $t\leq u$, with $u\geq2$, and
consider the case $t=u+1$. \ Thus $r+s=u+1\;(\leq n)$, and we may assume
$r,s\geq1$. \ When $r\geq2$,
\begin{align*}
\left\langle Y^{i}X^{j},Y^{r}X^{s}\right\rangle  &  =\left\langle Y^{i}%
X^{j},Y^{r-1}X^{s-1}\right\rangle \text{ \ (by Lemma \ref{lem1})}\\
&  =\left\langle Y^{i+1}X^{j-1},Y^{r-2}X^{s}\right\rangle \text{ \ (by the
inductive step)}\\
&  =\left\langle Y^{i+1}X^{j-1},Y^{r-1}X^{s+1}\right\rangle \text{ \ (by Lemma
\ref{lem1}),}%
\end{align*}
as desired. \ When $r=1$ and $j\geq2$, we have $s\leq n-r=n-1$, and we
consider three subcases.%

\noindent
\textbf{Subcase 1}. \ For $j=2$, $i=n-1$,
\begin{align*}
\left\langle Y^{n-1}X^{2},YX^{s}\right\rangle  &  =\left\langle Y^{n-2}%
X,YX^{s}\right\rangle \;\;\text{(by (\ref{intermediate}))}\\
&  =\left\langle Y^{n-2}X,X^{s-1}\right\rangle \;\;\text{(by Lemma
\ref{newlem}(iii))}\\
&  =\left\langle Y^{n-2},X^{s}\right\rangle \\
&  =\left\langle Y^{n-1},X^{s+1}\right\rangle \;\;\text{(by Lemma \ref{basic},
since }s\leq n-1\text{)}\\
&  =\left\langle Y^{n}X,X^{s+1}\right\rangle \;\;\text{(by (\ref{intermediate}%
)).}%
\end{align*}
\textbf{Subcase 2}. \ For $j\geq3,$ $i\geq1$,%
\begin{align*}
\left\langle Y^{i}X^{j},YX^{s}\right\rangle  &  =\left\langle Y^{i-1}%
X^{j-1},YX^{s}\right\rangle \text{ (by (\ref{intermediate}))}\\
&  =\left\langle Y^{i-1}X^{j-1},X^{s-1}\right\rangle \text{ (by Lemma
\ref{newlem}(iii))}\\
&  =\left\langle Y^{i-1}X^{j-3},X^{s+1}\right\rangle =\left\langle
Y^{i}X^{j-2},X^{s+1}\right\rangle \text{ (since }YX=1\text{ in }%
\mathcal{M}(n)\text{)}\\
&  =\left\langle Y^{i+1}X^{j-1},X^{s+1}\right\rangle \text{ (by
(\ref{intermediate})).}%
\end{align*}
\textbf{Subcase 3}. \ For $j=n+1,$ $i=0$,
\begin{align*}
\left\langle X^{n+1},YX^{s}\right\rangle  &  =\left\langle \frac{1}{a_{n}%
}[Y^{n-1}-Xp_{n-1}(X)-Xq_{n-1}(Y)],YX^{s}\right\rangle \\
&  =\left\langle \frac{1}{a_{n}}[Y^{n}-p_{n-1}(X)-q_{n-1}(Y)],X^{s}%
\right\rangle \text{ (by Lemma \ref{basic})}\\
&  =\left\langle X^{n},X^{s}\right\rangle =\left\langle X^{n-1},X^{s+1}%
\right\rangle \\
&  =\left\langle YX^{n},X^{s+1}\right\rangle \text{ (by (\ref{intermediate}%
)).}%
\end{align*}

Finally, when $r=1$ and $j=1$, we have $i=n\,$\ and $s\leq n-1$, so
\begin{align*}
\left\langle Y^{n}X,YX^{s}\right\rangle  &  =\left\langle Y^{n-1}%
,YX^{s}\right\rangle \;\;\text{(by (\ref{intermediate}))}\\
&  =\left\langle Y^{n},X^{s}\right\rangle =\left\langle a_{n}X^{n}%
+p_{n-1}(X)+q_{n-1}(Y),X^{s}\right\rangle \;\;\text{(by (\ref{yneq2}))}\\
&  =\left\langle a_{n}X^{n-1}+Yp_{n-1}(X)+Yq_{n-1}(Y),X^{s+1}\right\rangle
\;\;\\
&  \text{(by Lemma \ref{newlem}(ii) and Lemma \ref{basic}, as above)}\\
&  =\left\langle Y^{n+1},X^{s+1}\right\rangle .
\end{align*}
\end{proof}

Recall from (\ref{intermediate}) that columns $YX^{n},...,Y^{n}X$ are taken as
a block from consecutive columns of degree $n-1$ in $\mathcal{M}(n)$, so these
columns satisfy the Hankel property. \ Thus, in view of Lemma \ref{cor2}, the
next two results complete the proof that $B(n+1)$ has the Hankel property.

\begin{lemma}
\label{xn+1yk}For $k=1,...,n$,
\begin{equation}
\left\langle X^{n+1},Y^{k}\right\rangle =\left\langle YX^{n},Y^{k-1}%
X\right\rangle . \label{ynx}%
\end{equation}
\end{lemma}

\begin{proof}
We have
\begin{align*}
\left\langle X^{n+1},Y^{k}\right\rangle  &  =\left\langle \frac{1}{a_{n}%
}[Y^{n-1}-Xp_{n-1}(X)-Xq_{n-1}(Y)],Y^{k}\right\rangle \text{ \ (by
(\ref{xn+1}))}\\
&  =\left\langle \frac{1}{a_{n}}[Y^{n}-p_{n-1}(X)-q_{n-1}(Y)],Y^{k-1}%
\right\rangle \;\;\\
&  \text{(by Lemma \ref{newlem}(ii) and Lemma \ref{basic}, as above)}\\
&  =\left\langle X^{n},Y^{k-1}\right\rangle =\left\langle X^{n-1}%
,Y^{k-1}X\right\rangle \\
&  =\left\langle YX^{n},Y^{k-1}X\right\rangle \;\;\text{(by
(\ref{intermediate})).}%
\end{align*}
\end{proof}

\begin{lemma}
\label{yn+1xk}For $k=1,...,n$,
\[
\left\langle Y^{n}X,Y^{k}\right\rangle =\left\langle Y^{n+1},Y^{k-1}%
X\right\rangle .
\]
\end{lemma}

\begin{proof}
We have
\begin{align*}
\left\langle Y^{n}X,Y^{k}\right\rangle  &  =\left\langle Y^{n-1}%
,Y^{k}\right\rangle \;\;\text{(by (\ref{intermediate})}\\
&  =\left\langle Y^{n},Y^{k-1}\right\rangle =\left\langle a_{n}X^{n}%
+p_{n-1}(X)+q_{n-1}(Y),Y^{k-1}\right\rangle \text{ (by (\ref{yneq2}))}\\
&  =\left\langle a_{n}X^{n-1}+Yp_{n-1}(X)+Yq_{n-1}(Y),Y^{k-1}X\right\rangle
\text{ }\\
&  \text{(by Lemma \ref{basic} for the last two terms)}\\
&  =\left\langle Y^{n+1},Y^{k-1}X\right\rangle \;\;\text{(by (\ref{yn+1}))}%
\end{align*}
\end{proof}

The proof that block $B$ is of the form $\{B_{i,n+1}\}_{i=0}^{n}$ is now
complete. \ To finish the proof of Proposition \ref{propcaseIV} it now
suffices to show that $C:=W^{\ast}\mathcal{M}(n)W$ is Hankel. \ To do this,
observe that in the $C$ block of $\mathcal{M}:=\left[  \mathcal{M}%
(n);B\right]  =\left(
\begin{array}
[c]{cc}%
\mathcal{M}(n) & B\\
B^{\ast} & C
\end{array}
\right)  $, we need to compute inner products of the form $\left\langle
Y^{i}X^{j},Y^{k}X^{\ell}\right\rangle \;\;(i+j=k+\ell=n+1)$. \ For this, we
require an auxiliary lemma. \ For $i+j=k+\ell=n+1$, by $\left\langle
Y^{i}X^{j},Y^{k}X^{\ell}\right\rangle $ we mean, as usual, the entry in row
$Y^{k}X^{\ell}$ of column $Y^{i}X^{j}$; by self-adjointness of $\mathcal{M}%
(n+1)$, we have $\left\langle Y^{i}X^{j},Y^{k}X^{\ell}\right\rangle
=\left\langle Y^{k}X^{\ell},Y^{i}X^{j}\right\rangle $. \ Now suppose
$Y^{i}X^{j}=p(X,Y)$ and $Y^{k}X^{\ell}=q(X,Y)$, where $p(x,y)\equiv\sum_{0\leq
r+s\leq n}a_{rs}y^{r}x^{s}$ and $q(x,y)\equiv\sum_{0\leq t+u\leq n}b_{tu}%
y^{t}x^{u}$; we define $\left\langle p(X,Y),q(X,Y)\right\rangle :=$%
\newline $\sum_{0\leq r+s,t+u\leq n}a_{rs}b_{tu}\left\langle Y^{r}X^{s}%
,Y^{t}X^{u}\right\rangle $.

\begin{lemma}
\label{lem212}For $i+j,k+\ell=n+1$, $\left\langle Y^{i}X^{j},Y^{k}X^{\ell
}\right\rangle =\left\langle p(X,Y),q(X,Y)\right\rangle $.
\end{lemma}

\begin{proof}%
\begin{align*}
\left\langle p(X,Y),q(X,Y)\right\rangle  &  \equiv\sum_{0\leq r+s,t+u\leq
n}a_{rs}b_{tu}\left\langle Y^{r}X^{s},Y^{t}X^{u}\right\rangle \\
&  =\sum_{0\leq t+u\leq n}b_{tu}\left\langle \sum_{0\leq r+s\leq n}a_{rs}%
Y^{r}X^{s},Y^{t}X^{u}\right\rangle \\
&  =\sum_{0\leq t+u\leq n}b_{tu}\left\langle p(X,Y),Y^{t}X^{u}\right\rangle
=\sum_{0\leq t+u\leq n}b_{tu}\left\langle Y^{i}X^{j},Y^{t}X^{u}\right\rangle
\\
&  =\sum_{0\leq t+u\leq n}b_{tu}\left\langle Y^{t}X^{u},Y^{i}X^{j}%
\right\rangle \;\;\text{(since }\mathcal{M}=\mathcal{M}^{t}\text{)}\\
&  =\left\langle \sum_{0\leq t+u\leq n}b_{tu}Y^{t}X^{u},Y^{i}X^{j}%
\right\rangle =\left\langle q(X,Y),Y^{i}X^{j}\right\rangle \\
&  =\left\langle Y^{k}X^{\ell},Y^{i}X^{j}\right\rangle =\left\langle
Y^{i}X^{j},Y^{k}X^{\ell}\right\rangle \;\;\text{(since }\mathcal{M}%
=\mathcal{M}^{t}\text{).}%
\end{align*}
\end{proof}

\begin{proof}
[Proof of Proposition \textbf{\ref{propcaseIV}}]Note that since $M$ is a flat
extension, dependence relations in the columns of\ $(\mathcal{M}(n)\;\;B)$
extend to column relations in $(B^{\ast}\;\;C)$. \ In particular, the middle
$n$ columns of $C$ coincide with the columns of degree $n-1$ of $B^{\ast}$;
\ since $B$ has the Hankel property, so does $B^{\ast}$, and thus the middle
$n$ columns of $C$ have the Hankel property. \ To verify that $C$ is Hankel it
now suffices to focus on the first two and the last two columns of $C$, namely
$X^{n+1}$ and $YX^{n}$, and $Y^{n}X$ and $Y^{n+1}$. \ Since $C$ is
self-adjoint, and the middle $n$ columns have the Hankel property, to check
that $C$ is Hankel it only remains to show that $C_{n+2,1}=C_{n+1,2}$, i.e.,
$\left\langle X^{n+1},Y^{n+1}\right\rangle =\left\langle YX^{n},Y^{n}%
X\right\rangle $. \ Now, by (\ref{xn+1}), (\ref{yn+1}) and Lemma \ref{lem212}
we have
\begin{align*}
&  \left\langle X^{n+1},Y^{n+1}\right\rangle \\
&  =\left\langle \frac{1}{a_{n}}\{Y^{n-1}-X[p_{n-1}(X)+q_{n-1}(Y)]\},a_{n}%
X^{n-1}+Y[p_{n-1}(X)+q_{n-1}(Y)]\right\rangle \\
&  =\left\langle Y^{n-1},X^{n-1}\right\rangle +\frac{1}{a_{n}}\left\langle
Y^{n-1},Y[p_{n-1}(X)+q_{n-1}(Y)]\right\rangle \\
&  -\left\langle X[p_{n-1}(X)+q_{n-1}(Y)],X^{n-1}\right\rangle \\
&  -\frac{1}{a_{n}}\left\langle X[p_{n-1}(X)+q_{n-1}(Y)],Y[p_{n-1}%
(X)+q_{n-1}(Y)]\right\rangle \\
&  =\left\langle Y^{n},X^{n}\right\rangle +\frac{1}{a_{n}}\left\langle
Y^{n},p_{n-1}(X)+q_{n-1}(Y)\right\rangle -\left\langle p_{n-1}(X)+q_{n-1}%
(Y),X^{n}\right\rangle \\
&  -\frac{1}{a_{n}}\left\langle p_{n-1}(X)+q_{n-1}(Y)],p_{n-1}(X)+q_{n-1}%
(Y)]\right\rangle \\
&  \text{(by Lemma \ref{basic} for the first and fourth terms, and Lemma
\ref{newlem}(ii) for }\\
&  \text{the second and third terms)}%
\end{align*}%
\begin{align*}
&  =\left\langle \frac{1}{a_{n}}\{Y^{n}-[p_{n-1}(X)+q_{n-1}(Y)]\},a_{n}%
X^{n}+p_{n-1}(X)+q_{n-1}(Y)\right\rangle \\
&  =\left\langle X^{n},Y^{n}\right\rangle \text{ \ (by (\ref{xn}),
(\ref{yneq2}), and Lemma \ref{newlem}(i),(iii))}\\
&  =\left\langle X^{n-1},Y^{n-1}\right\rangle \text{ \ (by Lemma \ref{basic}%
)}\\
&  =\left\langle YX^{n},Y^{n}X\right\rangle \;\;\text{(by (\ref{intermediate})
and Lemma \ref{lem212}).}%
\end{align*}
This concludes the proof of Proposition \ref{propcaseIV}.
\end{proof}

\begin{remark}
\label{rem213}It is important for the sequel to note that in the proof of
Proposition \ref{propcaseIV}, the variety condition $\operatorname*{rank}%
\mathcal{M}(n)\leq\operatorname*{card}\mathcal{V}(\mathcal{M}(n))$ was used
only to show that $a_{n}\neq0$ in (\ref{yneq2}). \ Thus, if $\mathcal{M}(n)$
is positive, recursively generated, satisfies $YX=1$ in $\mathcal{C}%
_{\mathcal{M}(n)}$, and if the first dependence relation in $\mathcal{S}%
_{n}(n)$ is of the form (\ref{yneq2}) with $a_{n}\neq0$, then we may conclude
that $\mathcal{M}(n)$ has a flat extension $\mathcal{M}(n+1)$.
\end{remark}

\begin{proposition}
\label{propcaseV}(Case IV: $\operatorname*{rank}\mathcal{M}(n)=2n+1$) \ Assume
that $\mathcal{M}(n)$ is positive, recursively generated, and satisfies $YX=1
$. \ Assume also that $\mathcal{S}_{n}(n)$ is a basis for $\mathcal{C}%
_{\mathcal{M}(n)}$. \ Then $\mathcal{M}(n)$ admits a flat extension
$\mathcal{M}(n+1)$, or $\mathcal{M}(n)$ admits a positive, recursively
generated extension $\mathcal{M}(n+1)$, with $\operatorname*{rank}%
\mathcal{M}(n+1)=2n+2$, and $\mathcal{M}(n+1)$ admits a flat extension
$\mathcal{M}(n+2)$.
\end{proposition}

\begin{proof}
Since $YX=1$, and to guarantee that $\mathcal{M}(n+1)$ is recursively
generated, we define the middle $n$ columns of a proposed $B$ block for
$\mathcal{M}(n+1)$ as $[YX^{n}]_{m(n)}:=X^{n-1}\in\mathcal{C}_{\mathcal{M}%
(n)}$, $\left[  Y^{2}X^{n-1}\right]  _{m(n)}:=YX^{n-2}\in\mathcal{C}%
_{\mathcal{M}(n)}$, ... , $\left[  Y^{n}X\right]  _{m(n)}:=Y^{n-1}%
\in\mathcal{C}_{\mathcal{M}(n)}$. \ Moreover, if we wish to make $B_{n,n+1}$
Hankel, it is clear that all but the entry $\left\langle X^{n+1}%
,X^{n}\right\rangle $ in the column $\left[  X^{n+1}\right]  _{m(n)}$ must be
given in terms of entries in $\mathcal{M}(n)$, and that all but the entry
$\left\langle Y^{n+1},Y^{n}\right\rangle $ in $\left[  Y^{n+1}\right]
_{m(n)}$ must be given in terms of entries in $\mathcal{M}(n)$. \ To handle
the remaining entries we introduce two parameters $p$ and $q$; concretely, for
$i+j=0,...,n $,
\begin{equation}
\left\langle X^{n+1},Y^{i}X^{j}\right\rangle :=\left\{
\begin{array}
[c]{cc}%
\left\langle YX^{n},Y^{i-1}X^{j+1}\right\rangle  & (1\leq i\leq n)\\
\beta_{0,n+j+1} & (i=0,\;0\leq j\leq n-1)\\
p & (i=0,\;j=n)
\end{array}
\right.  , \label{eeq1}%
\end{equation}%
\begin{equation}
\left\langle Y^{n+1},Y^{i}X^{j}\right\rangle :=\left\{
\begin{array}
[c]{cc}%
\left\langle Y^{n}X,Y^{i+1}X^{j-1}\right\rangle  & (1\leq j\leq n)\\
\beta_{n+1+i,0} & (0\leq i\leq n-1,j=0)\\
q & (i=n,j=0)
\end{array}
\right.  . \label{eeq2}%
\end{equation}
A positive extension $\mathcal{M}(n+1)$ entails Ran $B\subseteq\;$Ran
$\mathcal{M}(n)$, so in particular we must show that $\left[  X^{n+1}\right]
_{m(n)}$, $\left[  Y^{n+1}\right]  _{m(n)}\in$ Ran $\mathcal{M}(n)$. \ To this
end, note that since $N:=[\mathcal{M}(n)]_{\mathcal{S}_{n}(n)}>0$, there exist
vectors $\mathbf{f},\mathbf{g}\in\mathbb{R}^{2n+1}$ such that $N\mathbf{f}%
=\left[  X^{n+1}\right]  _{\mathcal{S}_{n}(n)}$ and $N\mathbf{g}=\left[
Y^{n+1}\right]  _{\mathcal{S}_{n}(n)}$. \ Let $\mathbf{F},\mathbf{G\in
}\mathbb{R}^{m(n)}$ be given by
\[
\left\langle \mathbf{F},Y^{i}X^{j}\right\rangle :=\left\{
\begin{array}
[c]{c}%
\left\langle \mathbf{f},Y^{i}X^{j}\right\rangle \text{ \ if }Y^{i}X^{j}%
\in\mathcal{S}_{n}(n)\text{ }\\
0\text{ \ otherwise}%
\end{array}
\right.
\]
and
\[
\left\langle \mathbf{G},Y^{i}X^{j}\right\rangle :=\left\{
\begin{array}
[c]{c}%
\left\langle \mathbf{g},Y^{i}X^{j}\right\rangle \text{ \ if }Y^{i}X^{j}%
\in\mathcal{S}_{n}(n)\text{ }\\
0\text{ \ otherwise}%
\end{array}
\right.  .\
\]
We observe, for future reference, that since $\mathbf{f}=N^{-1}\left[
X^{n+1}\right]  _{\mathcal{S}_{n}(n)}$, $\mathbf{f}$ is linear in $p$ (and
independent of $q$), and so also is $\mathbf{F}$; similarly, $\mathbf{g}$ and
$\mathbf{G}$ are linear in $q$ and independent of $p$.

\textbf{Claim}. \ $\mathcal{M}(n)\mathbf{F}=\left[  X^{n+1}\right]  _{m(n)}$;
equivalently,
\[
\left\langle \mathcal{M}(n)\mathbf{F,}Y^{i}X^{j}\right\rangle =\left\langle
\left[  X^{n+1}\right]  _{m(n)},Y^{i}X^{j}\right\rangle
\]
for each $(i,j)\in I_{u}:=\{(i,j):i+j\leq n$ and ($(i=u\leq j\leq n)$ or
$(j=u\leq i\leq n)$)$\}$, $u=0,1,...,[\frac{n}{2}]$. \ Our proof of Claim 1 is
by induction on $u$. \ For $u=0$, we consider $Z\equiv Y^{i}$ or $Z\equiv
X^{j}$ in $\mathcal{S}_{n}(n)$, so
\begin{align*}
\left\langle \mathcal{M}(n)\mathbf{F,}Z\right\rangle  &  =\sum_{Y^{k}X^{\ell
}\in\mathcal{C}_{\mathcal{M}(n)}}\left\langle Y^{k}X^{\ell},Z\right\rangle
\left\langle \mathbf{F,}Y^{k}X^{\ell}\right\rangle \\
&  =\sum_{Y^{k}X^{\ell}\in\mathcal{S}_{n}(n)}\left\langle Y^{k}X^{\ell
},Z\right\rangle \left\langle \mathbf{f,}Y^{k}X^{\ell}\right\rangle
+\sum_{Y^{k}X^{\ell}\not \in \mathcal{S}_{n}(n)}\left\langle Y^{k}X^{\ell
},Z\right\rangle \cdot0\\
&  =\left\langle N\mathbf{f},Z\right\rangle =\left\langle \left[
X^{n+1}\right]  _{m(n)},Z\right\rangle ,
\end{align*}
as desired. \ We must now deal with rows of the form $Y^{i}X^{j}\;(i,j\geq1)
$; that is, we must prove that $\left\langle \mathcal{M}(n)\mathbf{F,}%
Y^{i}X^{j}\right\rangle =\left\langle \left[  X^{n+1}\right]  _{m(n)}%
,Y^{i}X^{j}\right\rangle $ for $i,j\geq1$ and $i+j\leq n$. \ Assume that the
Claim is true for $u=k\;(0\leq k\leq\lbrack\frac{n}{2}]-1)$, and consider
$(i,j)\in I_{k+1}$. \ We have
\begin{align*}
\left\langle \left[  X^{n+1}\right]  _{m(n)},Y^{i}X^{j}\right\rangle  &
\equiv\left\langle X^{n+1},Y^{i}X^{j}\right\rangle =\left\langle
YX^{n},Y^{i-1}X^{j+1}\right\rangle \ \;\text{(by (\ref{eeq1})}\\
&  =\left\langle X^{n-1},Y^{i-1}X^{j+1}\right\rangle .
\end{align*}
On the other hand,
\begin{align*}
\left\langle \mathcal{M}(n)\mathbf{F,}Y^{i}X^{j}\right\rangle  &
=\left\langle \mathcal{M}(n)\mathbf{F,}Y^{i-1}X^{j-1}\right\rangle \text{
\ (by Lemma \ref{newlem}(iii))}\\
&  =\left\langle \left[  X^{n+1}\right]  _{m(n)},Y^{i-1}X^{j-1}\right\rangle
\text{ (by the inductive step)}\\
&  =\left\langle X^{n+1},Y^{i-1}X^{j-1}\right\rangle .
\end{align*}
It thus suffices to prove that $\left\langle X^{n+1},Y^{i-1}X^{j-1}%
\right\rangle =\left\langle X^{n-1},Y^{i-1}X^{j+1}\right\rangle $. \ For
$i=1$,
\begin{align*}
\left\langle X^{n+1},X^{j-1}\right\rangle  &  =\beta_{0,n+j}\text{ (by
(\ref{eeq1}))}\\
&  =\left\langle X^{n-1},X^{j+1}\right\rangle ,
\end{align*}
and for $i>1$,
\begin{align*}
\left\langle X^{n+1},Y^{i-1}X^{j-1}\right\rangle  &  =\left\langle
YX^{n},Y^{i-2}X^{j}\right\rangle \text{ (by (\ref{eeq1})}\\
&  =\left\langle X^{n-1},Y^{i-2}X^{j}\right\rangle =\left\langle
X^{n-1},Y^{i-1}X^{j+1}\right\rangle \text{ (by Lemma \ref{newlem}(iii).}%
\end{align*}
This completes the proof of the Claim. \ An entirely similar argument, using
$\mathbf{g}$ instead of $\mathbf{f}$ and (\ref{eeq2}) instead of (\ref{eeq1}),
shows that $\mathcal{M}(n)\mathbf{G}=\left[  Y^{n+1}\right]  _{m(n)}$.
\ Moreover, by definition, $\left[  Y^{i}X^{j}\right]  _{m(n)}=Y^{i-1}%
X^{j-1}\in\mathcal{C}_{\mathcal{M}(n)}\;(i+j=n+1;$ $i,j\geq1)$, so we now have
Ran $B\subseteq$ Ran $\mathcal{M}(n)$; in particular, there exists $W$ such
that $\mathcal{M}(n)W=B$.

We note the following for future reference. \ From Lemma \ref{lem2.3} and the
fact that $\mathcal{M}(n)=\mathcal{M}(n)^{t}$, each row of $\mathcal{M}(n)$
coincides with a row indexed by an element of $\mathcal{S}_{n}(n)$. \ Since
$B=\mathcal{M}(n)W$, it now follows that each row of $(\mathcal{M}(n)\;B)$
coincides with a row of $(\mathcal{M}(n)\;B)$ indexed by an element of
$\mathcal{S}_{n}(n)$. \ 

We now form the flat extension $\mathcal{M}:=[\mathcal{M}(n);B]\equiv\left(
\begin{array}
[c]{cc}%
\mathcal{M}(n) & B\\
B^{t} & C
\end{array}
\right)  $, where $C:=W^{t}\mathcal{M}(n)W$. \ Exactly as in the proof of
Proposition \ref{propcaseIV}, $C$ is of the form
\[
C\equiv\left(
\begin{array}
[c]{ccccc}%
\tau & \beta_{0,2n} & \cdots & \beta_{02} & \eta\\
\beta_{0,2n} & \beta_{0,2n-2} & \cdots & \beta_{00} & \beta_{20}\\
\vdots & \vdots & \ddots & \vdots & \vdots\\
\beta_{02} & \beta_{00} & \cdots & \beta_{2n-2,0} & \beta_{2n,0}\\
\eta & \beta_{20} & \cdots & \beta_{2n,0} & \rho
\end{array}
\right)  ,
\]
where $C_{11}\equiv\tau:=\left[  X^{n+1}\right]  _{\mathcal{S}_{n}(n)}%
^{t}N^{-1}\left[  X^{n+1}\right]  _{\mathcal{S}_{n}(n)}$ and $C_{1,n+2}%
=C_{n+2,1}\equiv\eta:=\left[  X^{n+1}\right]  _{\mathcal{S}_{n}(n)}^{t}%
N^{-1}\left[  Y^{n+1}\right]  _{\mathcal{S}_{n}(n)}$. \ Thus, if $\eta
=\beta_{00}$, then $\mathcal{M}$ is a flat moment matrix extension of the form
$\mathcal{M}(n+1)$, and we are done.

Assume now that $\eta\neq\beta_{00}$. \ Let $u>\tau$ be arbitrary, and
consider the moment matrix $\mathcal{M}^{\prime}\equiv\mathcal{M}%
(n+1)^{\prime}$ obtained from $\mathcal{M}$ by replacing $\tau$ by $u$ and
$\eta$ by $\beta_{00}$. \ We partition $\mathcal{M}^{\prime}$ as\ $\mathcal{M}%
^{\prime}\equiv\left(
\begin{array}
[c]{cc}%
\mathcal{\tilde{M}} & \tilde{B}\\
\tilde{B}^{t} & \tilde{C}%
\end{array}
\right)  $, where $\mathcal{\tilde{M}}$ is the compression of $\mathcal{M}%
^{\prime}$ to rows and columns indexed by $\mathcal{\tilde{B}}:=\{1,X,Y,X^{2}%
,$ $YX,Y^{2},...,$ $X^{n},YX^{n-1}$,..., $Y^{n-1}X,Y^{n},X^{n+1}%
\}\subseteq\mathcal{C}_{\mathcal{M}^{\prime}}$ (i.e., $\mathcal{\tilde{M}}$ is
the extension of $\mathcal{M}(n)$ by row $X^{n+1}$ and column $X^{n+1}$ of
$\mathcal{M}^{\prime}$).

We claim that Ran $\tilde{B}\subseteq$ Ran $\mathcal{\tilde{M}}$. \ By the
flat construction of $\mathcal{M}$, the ``middle'' columns of $\left(
\begin{array}
[c]{c}%
B\\
C
\end{array}
\right)  $ are borrowed from columns of $\left(
\begin{array}
[c]{c}%
\mathcal{M}(n)\\
B^{t}%
\end{array}
\right)  $ of degree $n-1$, so, in particular, the columns of $\tilde{B}$
(except the rightmost column) are borrowed from columns in $\mathcal{\tilde
{M}}$. \ To prove the claim, it thus suffices to show that \newline $\left[
Y^{n+1}\right]  _{\mathcal{\tilde{B}}}\in$ Ran $\mathcal{\tilde{M}}$. \ Since
$u>\tau$, $\mathcal{\tilde{S}}:=\{1,X,Y,X^{2},Y^{2},...,X^{n},Y^{n},X^{n+1}\}$
is a basis for $\mathcal{C}_{\mathcal{\tilde{M}}}$, and $[\mathcal{\tilde{M}%
}]_{\mathcal{\tilde{S}}}$ is positive and invertible. \ Thus there exist
unique scalars $a_{1},a_{2},...,a_{2n+2}$ such that, in $\mathcal{C}%
_{(\mathcal{\tilde{M}}\;\tilde{B})}$, we have
\[
\lbrack Y^{n+1}]_{\mathcal{\tilde{S}}}=a_{1}[1]_{\mathcal{\tilde{S}}}%
+a_{2}[X]_{\mathcal{\tilde{S}}}+...+a_{2n+1}[Y^{n}]_{\mathcal{\tilde{S}}%
}+a_{2n+2}[X^{n+1}]_{\mathcal{\tilde{S}}}.\
\]
From the first part of the proof (concerning block $B$), we know that each row
of $(\mathcal{\tilde{M}}\;\tilde{B})$ coincides with a row indexed by an
element of $\mathcal{\tilde{S}}$, so it now follows that, in $\mathcal{C}%
_{(\mathcal{\tilde{M}}\;\tilde{B})}$,
\begin{equation}
\lbrack Y^{n+1}]_{\mathcal{\tilde{B}}}=a_{1}[1]_{\mathcal{\tilde{B}}}%
+a_{2}[X]_{\mathcal{\tilde{B}}}+...+a_{2n+1}[Y^{n}]_{\mathcal{\tilde{B}}%
}+a_{2n+2}[X^{n+1}]_{\mathcal{\tilde{B}}}, \label{star}%
\end{equation}
whence the claim is proved.

Since $\mathcal{\tilde{M}}\geq0$ and Ran $\tilde{B}\subseteq$ Ran
$\mathcal{\tilde{M}}$, we may construct the (positive) flat extension
$\mathcal{M}^{\flat}:=[\mathcal{\tilde{M}};\tilde{B}]\equiv\left(
\begin{array}
[c]{cc}%
\mathcal{\tilde{M}} & \tilde{B}\\
\tilde{B}^{t} & D
\end{array}
\right)  $, which we may re-partition as the moment matrix $\mathcal{M}%
(n+1)=\left(
\begin{array}
[c]{cc}%
\mathcal{M}(n) & B\\
B^{t} & C^{\flat}%
\end{array}
\right)  $, where $C^{\flat}$ is obtained from $C$ by replacing $\tau$ by $u
$, $\eta$ by $\beta_{00}$, and $\rho$ by some $\rho^{\flat}$ (determined by
extending (\ref{star}) to the full columns of $\mathcal{M}^{\flat}$).

Now $\mathcal{M}(n+1)$ is positive, recursively generated, satisfies $YX=1$,
and (by flatness of $\mathcal{M}^{\flat}$), $\operatorname*{rank}%
\mathcal{M}(n+1)=\operatorname*{rank}\mathcal{\tilde{M}}%
=1+\operatorname*{rank}\mathcal{M}(n)$. \ In $\mathcal{S}_{n+1}(n+1)$, the
first dependence relation is of the form $Y^{n+1}=a_{1}1+a_{2}X+...+a_{2n+1}%
Y^{n}+a_{2n+2}X^{n+1}$, and we assert that $a_{2n+2}\neq0$. \ Indeed, if
$a_{2n+2}=0$, then $[Y^{n+1}]_{\mathcal{S}_{n}(n)}=a_{1}[1]_{\mathcal{S}%
_{n}(n)}+a_{2}[X]_{\mathcal{S}_{n}(n)}+...+a_{2n+1}[Y^{n}]_{\mathcal{S}%
_{n}(n)}$, whence $(a_{1},...,a_{2n+1})^{t}=N^{-1}[Y^{n+1}]_{\mathcal{S}%
_{n}(n)}$. \ Now we have
\begin{align*}
\beta_{00}  &  =\left\langle [Y^{n+1}]_{\mathcal{\tilde{B}}},X^{n+1}%
\right\rangle \\
&  =a_{1}\left\langle [1]_{\mathcal{\tilde{B}}},X^{n+1}\right\rangle
+a_{2}\left\langle [X]_{\mathcal{\tilde{B}}},X^{n+1}\right\rangle
+...+a_{2n+1}\left\langle [Y^{n}]_{\mathcal{\tilde{B}}},X^{n+1}\right\rangle
\\
&  =[X^{n+1}]_{\mathcal{S}_{n}(n)}^{t}\cdot(a_{1},...,a_{2n+1})^{t}\\
&  =[X^{n+1}]_{\mathcal{S}_{n}(n)}^{t}N^{-1}[Y^{n+1}]_{\mathcal{S}_{n}%
(n)}=\eta,
\end{align*}
a contradiction. \ Since $a_{2n+2}\neq0$, we may now proceed exactly as in the
proof of Proposition \ref{propcaseIV} (beginning at (\ref{xn}) and replacing
$n$ by $n+1$) to conclude that $\mathcal{M}(n+1)$ admits a flat extension
$\mathcal{M}(n+2)$ (cf. Remark \ref{rem213}).
\end{proof}

\begin{remark}
Recall that $\mathbf{F}$ depends on $p$ and is independent of $q$, while
$\mathbf{G}$ depends on $q$ and is independent of $p$. \ It follows that
$\eta$ is of the form $\eta\equiv\eta(p,q)=a+bp+cq+dpq$, where $a,b,c,d\in
\mathbb{R}$ are independent of $p$ and $q$. \ Thus, if $b,$ $c$ or $d$ is
nonzero, it is possible to choose $p$ and $q$ so that $\eta=\beta_{00}$,
whence $\mathcal{M}(n)$ admits a flat extension $\mathcal{M}(n+1)$ (and
$\beta$ admits a $\operatorname*{rank}\mathcal{M}(n)$-atomic representing
measure). \ In \cite[Proof of Proposition 5.3]{tcmp5} we showed that this is
the case in the quartic moment problem ($n=2$), where we always have $b$ or
$d$ nonzero. \ For $n>2$, we do not know whether it is always the case that $b
$, $c$ or $d$ is nonzero.
\end{remark}

\begin{proof}
[Proof of Theorem \textbf{\ref{thmflat}}]Straightforward from Propositions
\ref{propcaseI}, \ref{propcaseIII}, \ref{propcaseIV} and \ref{propcaseV}.
\end{proof}

\section{\label{degenerate}The Truncated Moment Problem on Degenerate Hyperbolas}

In this section we prove Theorem \ref{firstmain}\ for degenerate hyperbolas.
\ By Proposition \ref{prop17}, it suffices to consider the case $yx=0$, and
the necessity of the conditions in Theorem \ref{firstmain} is clear from
Section \ref{Int}. \ We establish sufficiency in the following result.

\begin{theorem}
\label{secondflat}Let $\beta\equiv\beta^{(2n)}:\beta_{00},\beta_{01}%
,\beta_{10},...,\beta_{0,2n},...,\beta_{2n,0}$ be a family of real numbers,
$\beta_{00}>0$, and let $\mathcal{M}(n)$ be the associated moment matrix.
$\ $Assume that $\mathcal{M}(n)$ is positive, recursively generated, and
satisfies $YX=0$ and $\operatorname*{rank}\;\mathcal{M}(n)\leq
\operatorname*{card}\;\mathcal{V}(\beta)$. Then $\operatorname*{rank}%
\mathcal{M}(n)\leq2n+1$. \ If $\operatorname*{rank}\mathcal{M}(n)\leq2n$, then
$\mathcal{M}(n)$ admits a flat extension (so $\beta$ admits a
$\operatorname*{rank}\mathcal{M}(n)$-atomic representing measure supported in
$yx=0$). \ If $\operatorname*{rank}\mathcal{M}(n)=2n+1$, then $\mathcal{M}(n)$
admits a positive, recursively generated extension $\mathcal{M}(n+1)$,
satisfying $2n+1\leq\operatorname*{rank}\mathcal{M}(n+1)\leq2n+2$, and
$\mathcal{M}(n+1)$ admits a flat extension $\mathcal{M}(n+2)$ (so $\beta$
admits a representing measure $\mu$ supported in $yx=0$, with $2n+1\leq
\operatorname*{card}\operatorname*{supp}\mu\leq2n+2$).
\end{theorem}

By \cite[Theorem 2.1]{tcmp2}, we know that a positive, recursively generated
moment matrix $\mathcal{M}(n)$ admits flat extensions when $\left\{
\mathit{1},X,Y\right\}  $ is linearly dependent in $\mathcal{C}_{\mathcal{M}%
(n)}$; in the sequel, we therefore assume that $\left\{  \mathit{1}%
,X,Y\right\}  $ is linearly independent. \ We begin with an elementary lemma
based on recursiveness and $YX=0$. \ For $1\leq k\leq n$ let
\[
\mathcal{S}_{n}(k):=\{1,X,Y,X^{2},Y^{2},...,X^{k},Y^{k}\}\subseteq
\mathcal{C}_{\mathcal{M}(n)}.
\]

\begin{lemma}
For $n\geq2$, let $\mathcal{M}(n)$ be positive and recursively generated, and
assume that $YX=0$. \ Then each nonzero column of $\mathcal{M}(n)$ is in
$\mathcal{S}_{n}(n)$, and therefore $\operatorname*{rank}\mathcal{M}%
(n)\leq2n+1$.
\end{lemma}

In the sequel we also require the following well-known result.

\begin{lemma}
\label{inverse} (Choleski's Algorithm \cite{Atk}) Let $A$ be a positive and
invertible $d\times d$ matrix over $\mathbb{R}$, let $\mathbf{b}$ denote a
column vector in $\mathbb{R}^{d}$, and let $c\in\mathbb{R}$. \ Then
$\tilde{A}:=\left(
\begin{array}
[c]{cc}%
A & \mathbf{b}\\
\mathbf{b}^{t} & c
\end{array}
\right)  $ is positive and invertible if and only if $\delta:=c-\mathbf{b}%
^{t}A^{-1}\mathbf{b}>0$. \ In this case,
\[
\tilde{A}^{-1}=\frac{1}{\delta}\left(
\begin{array}
[c]{cc}%
(\delta\mathbf{+}A^{-1}\mathbf{bb}^{t})A^{-1} & -A^{-1}\mathbf{b}\\
-\mathbf{b}^{t}A^{-1} & 1
\end{array}
\right)  .
\]
\end{lemma}

We next divide the proof of Theorem \ref{secondflat} into three cases, based
on possible dependence relations among the elements of $\mathcal{S}_{n}(n)$.
\ Section \ref{examples} contains examples illustrating these cases. \ As in
Section \ref{3tm}, in each case, once we establish a flat extension, the
existence of the required representing measure $\mu$, necessarily supported in
$yx=0$, always follows immediately from (\ref{15}) and (\ref{exp18}), so we
will not repeat this argument in each case.

\begin{proposition}
\label{degcaseI}Suppose $\mathcal{M}(n)\equiv\mathcal{M}(n)(\beta)$ is
positive, recursively generated, $\operatorname*{card}\mathcal{V}%
(\mathcal{M}(n))\geq\operatorname*{rank}\mathcal{M}(n)$, and $YX=0$ in
$\mathcal{C}_{\mathcal{M}(n)}$. \ Suppose there exists $k$, $1<k\leq n$, such
that $\mathcal{S}_{n}(k-1)$ is linearly independent and $X^{k}\in$ lin.span
$\mathcal{S}_{n}(k-1)$. \ Then $\mathcal{M}(n)$ admits a flat extension
$\mathcal{M}(n+1) $ (and $\beta^{(2n)}$ admits a $\operatorname*{rank}%
\mathcal{M}(n)$-atomic representing measure).
\end{proposition}

\begin{proof}
By hypothesis, we may write
\begin{equation}
X^{k}=a_{0}1+a_{1}X+b_{1}Y+...+a_{k-1}X^{k-1}+b_{k-1}Y^{k-1}\;\;(a_{i}%
,b_{i}\in\mathbb{R)}. \label{degneweq1}%
\end{equation}
Equation (\ref{degneweq1}) implies that there are at most $k$ points in
$\mathcal{V}(\mathcal{M}(n))$ of the form $(x,0)$. \ If, for some $j$,
$b_{j}\neq0$, then it follows from (\ref{degneweq1}) that there are at most
$k-1$ points in $\mathcal{V}(\mathcal{M}(n))$ of the form $(0,y)$. \ In this
case, since $YX=0$ in $\mathcal{C}_{\mathcal{M}(n)}$, it follows that
$2k-1=\operatorname*{card}\mathcal{S}_{n}(k-1)\leq\operatorname*{rank}%
\mathcal{M}(n)\leq\operatorname*{card}\mathcal{V}(\mathcal{M}(n))\leq2k-1$,
whence $\operatorname*{rank}\mathcal{M}(n)=2k-1$ and $Y^{j}\in$ lin.span
$\mathcal{S}_{n}(k-1)$ ($k\leq j\leq n$). \ If $k=n$, this shows that
$\mathcal{M}(n)$ is flat. \ If $k<n$, then also $X^{j}\in$ lin.span
$\mathcal{S}_{n}(k-1)\;(k+1\leq j\leq n)$, so again $\mathcal{M}(n)$ is flat.
\ Thus, if some $b_{j}\neq0$, then $\mathcal{M}(n)$ is flat and the result
follows. \ We may thus assume that each $b_{j}=0$, i.e.,
\begin{equation}
X^{k}=a_{0}1+a_{1}X+...+a_{k-1}X^{k-1}. \label{degxk}%
\end{equation}
\ If $a_{0}\neq0$, there are no points in $\mathcal{V}(\mathcal{M}(n))$ of the
form $(0,y)$. \ Thus, in this case, each point in the variety is of the form
$(x,0)$, and (\ref{degxk}) implies that there can be at most $k$ such points.
\ Then $2k-1\leq\operatorname*{rank}\mathcal{M}(n)\leq\operatorname*{card}%
\mathcal{V}(\mathcal{M}(n))\leq k$, a contradiction. \ We thus conclude that
$a_{0}=0$, whence%
\begin{equation}
X^{k}=a_{1}X+...+a_{k-1}X^{k-1}. \label{degxk2}%
\end{equation}
By recursiveness,%
\begin{equation}
X^{k+i}=a_{1}X^{i+1}+...+a_{k-1}X^{k-1+i}\;\;(1\leq i\leq n-k). \label{degxk3}%
\end{equation}
Thus $\mathcal{\tilde{B}}:=\{1,X,Y,...,X^{k-1},Y^{k-1},Y^{k},...,Y^{n}\}$
spans $\mathcal{C}_{\mathcal{M}(n)}$. \ If $\mathcal{\tilde{B}}$ is linearly
dependent, it follows readily from recursiveness and $YX=0$ (as above) that
$\mathcal{M}(n)$ is flat, so again there is a (unique, $\operatorname*{rank}%
\mathcal{M}(n)$-atomic) representing measure.

We may thus assume that $\mathcal{\tilde{B}}$ is a basis for $\mathcal{C}%
_{\mathcal{M}(n)}$. \ We will show that $\mathcal{M}(n)$ admits infinitely
many flat extensions $\mathcal{M}(n+1)$ (each corresponding to a distinct
$\operatorname*{rank}\mathcal{M}(n)$-atomic representing measure). \ To
define
\[
\mathcal{M}(n+1)\equiv\left(
\begin{array}
[c]{cc}%
\mathcal{M}(n) & B(n+1)\\
B(n+1)^{t} & C(n+1)
\end{array}
\right)  ,
\]
which must be recursively generated (since $\mathcal{M}(n+1)$ would have a
representing measure by (\ref{exp18})), we first use (\ref{degxk2}) and
(\ref{degxk3}) to define $X^{n+1}$ in $\mathcal{C}_{(\mathcal{M}(n)\;B(n+1))}$
by
\begin{equation}
X^{n+1}:=a_{1}X^{n+2-k}+...+a_{k-1}X^{n}. \label{eq5}%
\end{equation}
Clearly, $X^{n+1}\in$ Ran $\mathcal{M}(n)$. \ To ensure the moment matrix
structure of $X^{n+1}$, we must verify that
\begin{equation}
\left\langle X^{n+1},Y^{i}X^{j}\right\rangle =\beta_{i,n+1+j}\;\;(0\leq
i+j\leq n-1). \label{eq6}%
\end{equation}
If $i>0$, then $\left\langle X^{p},Y^{i}X^{j}\right\rangle =0\;\;(0<p\leq n)$,
so (\ref{eq5}) implies $\left\langle X^{n+1},Y^{i}X^{j}\right\rangle
=0=\beta_{i,n+1+j}$ in this case. \ For $i=0$,
\begin{align*}
\left\langle X^{n+1},X^{j}\right\rangle  &  =\sum_{s=1}^{k-1}a_{s}\left\langle
X^{n+1-k+s},X^{j}\right\rangle =\sum_{s=1}^{k-1}a_{s}\left\langle
X^{n-k+s},X^{j+1}\right\rangle \\
&  =\left\langle X^{n},X^{j+1}\right\rangle \;\;\text{(by (\ref{degxk3}))}\\
&  =\beta_{0,n+1+j}\;\;\text{(since }0\leq j\leq n-1\text{).}%
\end{align*}

We next define $Y^{n+1}$ in $\mathcal{C}_{(\mathcal{M}(n)\;B(n+1))}$
consistent with moment matrix structure and $YX=0$. \ Let $\mathcal{B}%
:=\{1,X,...,X^{k-1}\}\bigcup\{Y,...,,Y^{k-1},...,Y^{n-1}\}\subseteq
\mathcal{C}_{\mathcal{M}(n)}$; $\left[  Y^{n+1}\right]  _{\mathcal{B}}$ is
uniquely determined by $\left\langle Y^{n+1},X^{j}\right\rangle =0\;(1\leq
j\leq k-1)$ and $\left\langle Y^{n+1},Y^{i}\right\rangle =\beta_{n+1+i,0}%
\;(0\leq i\leq n-1)$. \ Now
\begin{equation}
\left[  Y^{n+1}\right]  _{\mathcal{\tilde{B}}}=\left(
\begin{array}
[c]{c}%
\left[  Y^{n+1}\right]  _{\mathcal{B}}\\
q
\end{array}
\right)  \;(q\in\mathbb{R)}, \label{eqnew65}%
\end{equation}
and every other component of $Y^{n+1}$ (besides those in $\left[
Y^{n+1}\right]  _{\mathcal{\tilde{B}}}$) must be zero (corresponding to
$\left\langle Y^{n+1},Y^{r}X^{s}\right\rangle $ for some $s>0$).

We will show that for each $q\in\mathbb{R}$, $Y^{n+1}\in$Ran $\mathcal{M}(n)$.
\ Let $\mathcal{M}:=[\mathcal{M}(n)]_{\mathcal{B}}>0$. \ For each vector of
the form $[v]_{\mathcal{B}}$, $\mathcal{M}^{-1}([v]_{\mathcal{B}})\equiv
(r_{0},r_{1},s_{1},...,r_{k-1},s_{k-1},...,s_{n-1})$ uniquely satisfies
\begin{equation}
\lbrack v]_{\mathcal{B}}=r_{0}[1]_{\mathcal{B}}+r_{1}[X]_{\mathcal{B}}%
+s_{1}[Y]_{\mathcal{B}}+...+r_{k-1}[X^{k-1}]_{\mathcal{B}}+s_{k-1}%
[Y^{k-1}]_{\mathcal{B}}+...+s_{n-1}[Y^{n-1}]_{\mathcal{B}}.\ \label{eqnew7}%
\end{equation}
We claim that
\begin{equation}
\left[  X^{k+j}\right]  _{\mathcal{B}}^{t}\mathcal{M}^{-1}\left[
Y^{n}\right]  _{\mathcal{B}}=0\;(0\leq j\leq n-k).\ \label{eqnew8}%
\end{equation}
Indeed, (\ref{degxk2}) and (\ref{degxk3}) together imply that $X^{k+j}$ is a
linear combination of $X$, $X^{2}$, ... , $X^{k-1}$, so in (\ref{eqnew7}), all
the coordinates of $\mathcal{M}^{-1}[X^{k+j}]_{\mathcal{B}}$ corresponding to
$1,Y,...,Y^{n-1}$ equal zero. \ By the $YX=0$ structure of $\mathcal{M}(n)$,
the only coordinates of $[Y^{n}]_{\mathcal{B}}$ that are possibly nonzero
correspond to $1,Y,...,Y^{n-1}$. \ These observations show that $\left[
Y^{n}\right]  _{\mathcal{B}}^{t}\mathcal{M}^{-1}\left[  X^{k+j}\right]
_{\mathcal{B}}=0$, whence (\ref{eqnew8}) follows. \ Similarly, the $YX=0$
structure of $Y^{n+1}$ implies that
\begin{equation}
\left[  X^{k+j}\right]  _{\mathcal{B}}^{t}\mathcal{M}^{-1}\left[
Y^{n+1}\right]  _{\mathcal{B}}=0\;(0\leq j\leq n-k). \label{eqnew9}%
\end{equation}
Since $\mathcal{\tilde{M}}:=[\mathcal{M}(n)]_{\mathcal{\tilde{B}}}>0$, we have
$\left[  Y^{n+1}\right]  _{\mathcal{\tilde{B}}}\in$ Ran $\mathcal{\tilde{M}}$.
\ To show that $Y^{n+1}\in$ Ran $\mathcal{M}(n)$, the $YX=0$ structure of
$\mathcal{M}(n)$ and of $Y^{n+1}$ imply that it suffices to verify that
\begin{equation}
\left[  X^{k+j}\right]  _{\mathcal{\tilde{B}}}^{t}(\mathcal{\tilde{M}}%
)^{-1}\left[  Y^{n+1}\right]  _{\mathcal{\tilde{B}}}=0\;(0\leq j\leq
n-k).\ \label{eqnew10}%
\end{equation}
Now $\mathcal{\tilde{M}}=\left(
\begin{array}
[c]{cc}%
\mathcal{M} & \left[  Y^{n}\right]  _{\mathcal{B}}\\
\left[  Y^{n}\right]  _{\mathcal{B}}^{t} & c
\end{array}
\right)  \equiv\left(
\begin{array}
[c]{cc}%
\mathcal{M} & \mathbf{b}\\
\mathbf{b}^{t} & c
\end{array}
\right)  $, so\
\[
(\mathcal{\tilde{M}})^{-1}=\frac{1}{\delta}\left(
\begin{array}
[c]{cc}%
(\delta+\mathcal{M}^{-1}\mathbf{bb}^{t})\mathcal{M}^{-1} & -\mathcal{M}%
^{-1}\mathbf{b}\\
-\mathbf{b}^{t}\mathcal{M}^{-1} & 1
\end{array}
\right)  ,
\]
with $\delta:=c-\mathbf{b}^{t}\mathcal{M}^{-1}\mathbf{b}$, by Lemma
\ref{inverse}. \ We also have $\left[  X^{k+j}\right]  _{\mathcal{\tilde{B}}%
}^{t}=\left(  \left[  X^{k+j}\right]  _{\mathcal{B}}^{t}\;\;0\right)  $ (since
$\left\langle X^{k+j},Y^{n}\right\rangle =0$); thus, using (\ref{eqnew65}), we
have
\begin{align*}
\left[  X^{k+j}\right]  _{\mathcal{\tilde{B}}}^{t}\mathcal{\tilde{M}}%
^{-1}\left[  Y^{n+1}\right]  _{\mathcal{\tilde{B}}}  &  =\frac{1}{\delta
}\left[  X^{k+j}\right]  _{\mathcal{B}}^{t}((\delta+\mathcal{M}^{-1}%
\mathbf{bb}^{t})\mathcal{M}^{-1}\left[  Y^{n+1}\right]  _{\mathcal{B}%
}-q\mathcal{M}^{-1}\mathbf{b})\\
&  =\left[  X^{k+j}\right]  _{\mathcal{B}}^{t}\mathcal{M}^{-1}\left[
Y^{n+1}\right]  _{\mathcal{B}}\\
&  +\frac{1}{\delta}(\left[  X^{k+j}\right]  _{\mathcal{B}}^{t}\mathcal{M}%
^{-1}\mathbf{b)b}^{t}\mathcal{M}^{-1}\left[  Y^{n+1}\right]  _{\mathcal{B}}\\
&  -\frac{q}{\delta}\left[  X^{k+j}\right]  _{\mathcal{B}}^{t}\mathcal{M}%
^{-1}\mathbf{b}\\
&  =0\;\;\text{(by (\ref{eqnew8}) and (\ref{eqnew9}).}%
\end{align*}
\ Thus $Y^{n+1}\in$ Ran $\mathcal{M}(n)$. \ 

We now have columns $X^{n+1}$ and $Y^{n+1}$ for block $B(n+1)$, and to
preserve the $YX=0$ structure we set $Y^{i}X^{j}=0\;(i+j=n+1;i,j>0)$. \ Thus
Ran $B(n+1)\subseteq$ Ran $\mathcal{M}(n)$, and we consider the flat
extension
\begin{equation}
\mathcal{M}^{\symbol{94}}:=\left[  \mathcal{M}(n);B(n+1)\right]
\equiv\left(
\begin{array}
[c]{cc}%
\mathcal{M}(n) & B(n+1)\\
B(n+1)^{t} & C
\end{array}
\right)  . \label{eqnew11}%
\end{equation}
To complete the proof, we will show that $C$ has the form of a moment matrix
block consistent with $YX=0$. \ To do so, from (\ref{eqnew11}) and the $YX=0$
structures of $\mathcal{M}(n)$ and $B(n+1)$, it suffices to show that
$\left\langle X^{n+1},Y^{n+1}\right\rangle =0$. \ Now,
\begin{align*}
\left\langle X^{n+1},Y^{n+1}\right\rangle  &  =\sum_{j=1}^{k-1}a_{j}%
\left\langle X^{n+1-k+j},Y^{n+1}\right\rangle \;\;\text{(by (\ref{eq5}) and
(\ref{eqnew11}))}\\
&  =\sum_{j=1}^{k-1}a_{j}\left\langle Y^{n+1},X^{n+1-k+j}\right\rangle
\;\;\text{(since }\mathcal{M}^{\symbol{94}}=(\mathcal{M}^{\symbol{94}}%
)^{t}\text{)}\\
&  =0\;\;\text{(since }Y^{n+1}\text{ in }B(n+1)\text{ has }YX=0\text{
structure).}%
\end{align*}
Thus $\mathcal{M}^{\symbol{94}}$ is a flat moment matrix extension of
$\mathcal{M}(n)$; the proof is complete.
\end{proof}

We next consider the case when the first column dependence relation occurs at
$Y^{k}\;(1<k\leq n)$. \ 

\begin{proposition}
\label{degcaseII}Assume $\mathcal{M}(n)(\beta)$ is positive, recursively
generated, and $YX=0$ in $\mathcal{C}_{\mathcal{M}(n)}$. \ Suppose that for
some $k\;(1<k\leq n)$, $\mathcal{S}\equiv\{1,X,Y,...,X^{k-1},Y^{k-1},X^{k}\} $
is linearly independent and $Y^{k}\in$ lin.span $\mathcal{S}$. \ Then
$\mathcal{M}(n)$ admits a flat extension (and $\beta^{(2n)}$ admits a
$\operatorname*{rank}\mathcal{M}(n)$-atomic representing measure).
\end{proposition}

\begin{proof}
In $\mathcal{C}_{\mathcal{M}(n)}$ we have a dependence relation
\begin{equation}
Y^{k}=a_{0}1+a_{1}X+b_{1}Y+...+a_{k-1}X^{k-1}+b_{k-1}Y^{k-1}+a_{k}X^{k}.
\label{eqq1}%
\end{equation}
\ If $a_{k}=0$, we can interchange the roles of $X$ and $Y$ and invoke
Proposition \ref{degcaseI}. \ We may thus assume $a_{k}\neq0$, so there is a
dependence relation of the form
\begin{equation}
X^{k}=\tilde{a}_{0}1+\tilde{a}_{1}X+\tilde{b}_{1}Y+...+\tilde{a}_{k-1}%
X^{k-1}+\tilde{b}_{k-1}Y^{k-1}+\tilde{b}_{k}Y^{k}. \label{eqq2}%
\end{equation}
Since any flat extension $\mathcal{M}(n+1)$ must be recursively generated,
with $YX=0$, in $B(n+1)$ we must have
\begin{equation}
Y^{n+1}=a_{0}Y^{n+1-k}+b_{1}Y^{n+2-k}+...+b_{k-1}Y^{n} \label{eqq3}%
\end{equation}
and
\begin{equation}
X^{n+1}=\tilde{a}_{0}X^{n+1-k}+\tilde{a}_{1}X^{n+2-k}+...+\tilde{a}_{k-1}%
X^{n}\text{;} \label{eqq4}%
\end{equation}
let $\mathbf{w}$ denote the column vector of length $m(n)$ such that
\[
\left\langle \mathbf{w},Y^{i}X^{j}\right\rangle :=\left\{
\begin{array}
[c]{l}%
\tilde{a}_{j-n-1+k}\;\;\text{if }n+1-k\leq j\leq n\text{ and }i=0\\
0\;\;\text{otherwise}%
\end{array}
.\right.
\]
To show that $Y^{n+1}$ is consistent with moment matrix structure and $YX=0$,
we first verify that
\begin{equation}
\left\langle Y^{n+1},Y^{i}X^{j}\right\rangle =\beta_{n+1+i,j}\;(0\leq i+j\leq
n-1)\text{;} \label{eqq5}%
\end{equation}
indeed,
\begin{align*}
\left\langle Y^{n+1},Y^{i}X^{j}\right\rangle  &  =\left\langle a_{0}%
Y^{n+1-k},Y^{i}X^{j}\right\rangle +\left\langle b_{1}Y^{n+2-k},Y^{i}%
X^{j}\right\rangle +...+\left\langle b_{k-1}Y^{n},Y^{i}X^{j}\right\rangle \\
&  =\left\langle a_{0}Y^{n-k},Y^{i+1}X^{j}\right\rangle +\left\langle
b_{1}Y^{n+1-k},Y^{i+1}X^{j}\right\rangle +...\\
&  +\left\langle b_{k-1}Y^{n-1},Y^{i+1}X^{j}\right\rangle \\
&  =\left\langle Y^{n},Y^{i+1}X^{j}\right\rangle \;\;\text{(by recursiveness,
using (\ref{eqq1})}\\
&  =\beta_{n+i+1,j}\;\;\text{(by the structure of }\mathcal{M}(n)\text{).}%
\end{align*}
Since (\ref{eqq3}) readily implies that $\left\langle Y^{n+1},Y^{i}%
X^{j}\right\rangle =0$ when $j>0$ and $i+j=n$, it follows that $Y^{n+1}$ is
consistent. \ A similar argument (using (\ref{eqq2}) and (\ref{eqq4})) shows
that $X^{n+1}$ is also consistent with moment matrix structure and $YX=0$.
\ Now, setting $Y^{i}X^{j}=0$ in $B(n+1)$ for $i+j=n+1\;(i,j>0)$, we have a
moment matrix block $B(n+1)$ consistent with $YX=0$ and satisfying Ran
$B(n+1)\subseteq$ Ran $\mathcal{M}(n)$. \ Consider $\mathcal{M}\symbol{94}%
:=\left[  \mathcal{M}(n);B(n+1)\right]  \equiv\left(
\begin{array}
[c]{cc}%
\mathcal{M}(n) & B(n+1)\\
B(n+1)^{t} & C
\end{array}
\right)  $. \ To show that $\mathcal{M}\symbol{94}$ is a moment matrix, it
suffices to check that in block $C$, $C_{n+2,1}\equiv\left\langle
X^{n+1},Y^{n+1}\right\rangle =0$, and this follows immediately from the
identity $C_{n+2,1}=[Y^{n+1}]_{m(n)}^{t}\cdot\mathbf{w}$.
\end{proof}

The following result concludes the proof of Theorem \ref{secondflat}.

\begin{proposition}
\label{degcaseIII}Assume that $\mathcal{M}\equiv\mathcal{M}(n)$ is positive,
recursively generated, satisfies $YX=0$, and that $\mathcal{S}_{n}(n)$ is a
basis for $\mathcal{C}_{\mathcal{M}(n)}$. \ Then either $\mathcal{M}(n)$
admits a flat extension (and $\beta$ admits a $(2n+1)$-atomic representing
measure) or $\mathcal{M}(n)$ admits a rank-$(2n+2)$ positive, recursively
generated extension $\mathcal{M}(n+1)$ which has a flat extension
$\mathcal{M}(n+2)$ (and $\beta$ admits a $(2n+2)$-atomic representing measure).
\end{proposition}

\begin{proof}
By hypothesis, $\mathcal{B}:=\{1,X,Y,...,X^{n},Y^{n}\}$ is a basis for
$\mathcal{S}_{n}(n)$, so the compression $A\equiv\left[  \mathcal{M}\right]
_{\mathcal{B}}$ of $\mathcal{M}$ to the rows and columns of $\mathcal{B}$ is
positive and invertible. \ We define the columns $X^{n+1}$ and $Y^{n+1}$ in
the proposed block $B(n+1)$ by
\[
\left\langle X^{n+1},X^{i}\right\rangle :=\left\{
\begin{array}
[c]{cc}%
\left\langle X^{n},X^{i+1}\right\rangle \;(=\beta_{0,n+i+1}) & 0\leq i<n\\
p & i=n
\end{array}
\right.
\]%
\[
\left\langle Y^{n+1},Y^{j}\right\rangle :=\left\{
\begin{array}
[c]{cc}%
\left\langle Y^{n},Y^{j+1}\right\rangle \;(=\beta_{n+j+1,0}) & 0\leq j<n\\
q & j=n
\end{array}
\right.
\]
and $\left\langle X^{n+1},Y^{j}X^{i}\right\rangle :=0\;(1\leq i+j\leq
n;j\geq1)$ and $\left\langle Y^{n+1},Y^{j}X^{i}\right\rangle :=0\;(1\leq
i+j\leq n;i\geq1)$, where $p$ and $q$ are two parameters. \ Let $\mathbf{r}%
_{p}:=\left[  X^{n+1}\right]  _{\mathcal{B}}$ and $\mathbf{s}_{q}:=\left[
Y^{n+1}\right]  _{\mathcal{B}}$. \ Due to the $YX=0$ structure of
$\mathcal{M}(n)$, it is straightforward to check that Ran $B(n+1)\subseteq$
Ran $\mathcal{M}$, and that if $W$ is a matrix satisfying $\mathcal{M}%
W=B(n+1)$, then
\[
C\equiv W^{\ast}\mathcal{M}W=\left(
\begin{array}
[c]{ccccc}%
\mathbf{r}_{p}^{t}A^{-1}\mathbf{r}_{p} & 0 & \cdots & 0 & \mathbf{r}_{p}%
^{t}A^{-1}\mathbf{s}_{q}\\
0 & 0 & \cdots & 0 & 0\\
\vdots & \vdots & \ddots & \vdots & \vdots\\
0 & 0 & \cdots & 0 & 0\\
\mathbf{s}_{q}^{t}A^{-1}\mathbf{r}_{p} & 0 & \cdots & 0 & \mathbf{s}_{q}%
^{t}A^{-1}\mathbf{s}_{q}%
\end{array}
\right)  .
\]
It follows at once that $\mathcal{M}$ admits a flat extension $\mathcal{M}%
(n+1)$ if and only if there exist real numbers $p$ and $q$ such that
\[
\alpha(p,q):=\mathbf{s}_{q}^{t}A^{-1}\mathbf{r}_{p}=0.\
\]

We may thus assume that $\alpha$ is nonzero on $\mathbb{R}^{2}$. \ Fix
$p,q\in\mathbb{R}$ and let $u>\mathbf{r}_{p}^{t}A^{-1}\mathbf{r}_{p}$ and
$v:=\frac{\alpha^{2}}{u-\mathbf{r}_{p}^{t}A^{-1}\mathbf{r}_{p}}+\mathbf{s}%
_{q}^{t}A^{-1}\mathbf{s}_{q}$, so that
\[
C(u,v):=\left(
\begin{array}
[c]{ccccc}%
u & 0 & \cdots & 0 & 0\\
0 & 0 & \cdots & 0 & 0\\
\vdots & \vdots & \ddots & \vdots & \vdots\\
0 & 0 & \cdots & 0 & 0\\
0 & 0 & \cdots & 0 & v
\end{array}
\right)  \geq C(n+1)
\]
and
\[
\operatorname*{rank}(C(u,v)-C)=1.
\]
Then
\[
\mathcal{M}(n+1)\equiv\mathcal{M}(n+1;u,p,q):=\left(
\begin{array}
[c]{cc}%
\mathcal{M} & B(n+1)\\
B(n+1)^{\ast} & C(u,v)
\end{array}
\right)
\]
is positive and recursively generated, and $\operatorname*{rank}%
\mathcal{M}(n+1)=1+\operatorname*{rank}\mathcal{M}$ ($=2n+2$).

We claim that $\mathcal{M}(n+1)$ admits a flat extension $\mathcal{M}(n+2)$.
\ We first show that there is a unique block $B(n+2)$, subordinate to $YX=0$,
such that Ran $B(n+2)\subseteq$ Ran $\mathcal{M}(n+1)$. \ In any such block,
$YX^{n+1}=Y^{2}X^{n}=...=Y^{n+1}X=0$. \ In $Y^{n+2}$, all the entries are
determined from $\mathcal{M}(n+1)$ and $YX=0$, except $s:=\left\langle
Y^{n+2},Y^{n+1}\right\rangle $. \ Now, in $\mathcal{C}_{\mathcal{M}(n+1)}$ we
have a dependence relation
\[
Y^{n+1}=p_{n+1}(X)+q_{n}(Y),
\]
with $\deg p_{n+1}\leq n+1,\deg q_{n}\leq n$. \ Since a flat extension must
necessarily be recursively generated, in $\mathcal{C}_{\mathcal{M}(n+2)}$ we
must have
\begin{equation}
Y^{n+2}=p_{n+1}(0)Y+Yq_{n}(Y), \label{degyn+2}%
\end{equation}
whence $s=\left\langle p_{n+1}(0)Y+Yq_{n}(Y),Y^{n+1}\right\rangle $. \ Thus,
$Y^{n+2}\in$ Ran $\mathcal{M}(n+1)$ and $Y^{n+2}$ has $YX=0$ structure and is
Hankel with respect to $Y^{n+1}X$.

Note that if $X^{n+2}$ for block $B(n+2)$ is defined to be consistent with
known moment values and $YX=0$ structure, then every value in $X^{n+2}$ is
determined except $r:=\left\langle X^{n+2},X^{n+1}\right\rangle $. \ We next
show that there is a unique value of $r$ such that $X^{n+2}\in$ Ran
$\mathcal{M}(n+1)$. \ Since $u>\mathbf{r}_{p}^{t}A^{-1}\mathbf{r}_{p}^{t}$,
Smul'jan's Theorem \cite{Smu} (described in Section \ref{Int}) implies that
the compression $\tilde{A}$ of $\mathcal{M}(n+1)$ to rows and columns given by
$\mathcal{\tilde{B}}:=\mathcal{B}\bigcup\{X^{n+1}\}$ is positive and
invertible, and of the form $\tilde{A}=\left(
\begin{array}
[c]{cc}%
A & \mathbf{r}_{p}\\
\mathbf{r}_{p}^{t} & u
\end{array}
\right)  $. \ We now apply Lemma \ref{inverse} to obtain
\[
\tilde{A}^{-1}=\frac{1}{\delta}\left(
\begin{array}
[c]{cc}%
(\delta\mathbf{+}A^{-1}\mathbf{r}_{p}\mathbf{r}_{p}^{t})A^{-1} &
-A^{-1}\mathbf{r}_{p}\\
-\mathbf{r}_{p}^{t}A^{-1} & 1
\end{array}
\right)  \;\;(\delta:=u-\mathbf{r}_{p}^{t}A^{-1}\mathbf{r}_{p}).
\]
Observe that $\left[  X^{n+2}\right]  _{\mathcal{\tilde{B}}}\equiv\left(
\begin{array}
[c]{c}%
\left[  X^{n+2}\right]  _{\mathcal{B}}\\
r
\end{array}
\right)  $, and that, apart from the entries in $\left[  X^{n+2}\right]
_{\mathcal{\tilde{B}}}$, all other entries of $X^{n+2}$ in $B(n+2)$ are zero.
\ From the $YX=0$ structure of $\mathcal{M}(n)$, it follows that $X^{n+2}\in$
Ran $\mathcal{M}(n+1)$ if and only if
\begin{equation}
\left\langle X^{n+2},Y^{n+1}\right\rangle :=\left[  Y^{n+1}\right]
_{\mathcal{\tilde{B}}}^{t}\tilde{A}^{-1}\left[  X^{n+2}\right]
_{\mathcal{\tilde{B}}}=0, \label{degeq1}%
\end{equation}
so it suffices to show that (\ref{degeq1}) admits a unique solution for $r$.
\ Now,
\begin{multline*}
\left[  Y^{n+1}\right]  _{\mathcal{\tilde{B}}}^{t}\tilde{A}^{-1}\left[
X^{n+2}\right]  _{\mathcal{\tilde{B}}}=\left(
\begin{array}
[c]{cc}%
\mathbf{s}_{q}^{t} & 0
\end{array}
\right)  \left(
\begin{array}
[c]{cc}%
A & \mathbf{r}_{p}\\
\mathbf{r}_{p}^{t} & u
\end{array}
\right)  ^{-1}\left(
\begin{array}
[c]{c}%
\left[  X^{n+2}\right]  _{\mathcal{B}}\\
r
\end{array}
\right) \\
=\left(
\begin{array}
[c]{cc}%
\mathbf{s}_{q}^{t} & 0
\end{array}
\right)  \frac{1}{\delta}\left(
\begin{array}
[c]{cc}%
(\delta\mathbf{+}A^{-1}\mathbf{r}_{p}\mathbf{r}_{p}^{t})A^{-1} &
-A^{-1}\mathbf{r}_{p}\\
-\mathbf{r}_{p}^{t}A^{-1} & 1
\end{array}
\right)  \left(
\begin{array}
[c]{c}%
\left[  X^{n+2}\right]  _{\mathcal{B}}\\
r
\end{array}
\right) \\
=\frac{1}{\delta}\left(
\begin{array}
[c]{cc}%
\mathbf{s}_{q}^{t} & 0
\end{array}
\right)  \left(
\begin{array}
[c]{c}%
H-rA^{-1}\mathbf{r}_{p}\\
\ast
\end{array}
\right)  \;\;\text{(for a certain vector }H\text{)}\\
=\frac{1}{\delta}(\mathbf{s}_{q}^{t}H-r\mathbf{s}_{q}^{t}A^{-1}\mathbf{r}%
_{p})=\frac{1}{\delta}(\mathbf{s}_{q}^{t}H-r\alpha(p,q))\text{.}%
\end{multline*}
Since $\alpha(p,q)\neq0$, it follows that (\ref{degeq1}) admits a unique
solution $r\equiv r(p,q,u)$. \ 

With this value, Ran $B(n+1)\subseteq$ Ran $\mathcal{M}(n+1)$, so
$B(n+2)=\mathcal{M}(n+1)W$ for some matrix $W$. \ To show that the flat
extension $\mathcal{M}_{n+2}:=\left[  \mathcal{M}(n+1);B(n+2)\right]
\equiv\left(
\begin{array}
[c]{cc}%
\mathcal{M}(n+1) & B(n+2)\\
B(n+2)^{t} & C_{n+2}%
\end{array}
\right)  $ is a moment matrix, it now suffices to show that \newline
$\left\langle Y^{n+2},X^{n+2}\right\rangle =0$; this is because, by positivity
of $\mathcal{M}_{n+2}$, $\left\langle X^{n+2},Y^{n+2}\right\rangle
=\left\langle Y^{n+2},X^{n+2}\right\rangle $ and, by flatness, all other
entries of $C_{n+2}$ (except $\left\langle X^{n+2},X^{n+2}\right\rangle $ and
$\left\langle Y^{n+2},Y^{n+2}\right\rangle $) are clearly zero. \ Now recall
that $Y^{n+2}=p_{n+1}(0)Y+Yq_{n}(Y)$ in $\mathcal{C}_{(\mathcal{M}%
(n+1)\;B(n+2))}$ (by (\ref{degyn+2})), so by flatness, the same relation must
hold in $\mathcal{C}_{\mathcal{M}_{n+2}}$. \ Thus,
\[
\left\langle Y^{n+2},X^{n+2}\right\rangle =\left\langle p_{n+1}(0)Y+Yq_{n}%
(Y),X^{n+2}\right\rangle =\left\langle X^{n+2},p_{n+1}(0)Y+Yq_{n}%
(Y)\right\rangle =0,
\]
since, by the construction of $X^{n+2}$ in $B(n+2)$ (consistent with the
relation $YX=0$), $\left\langle X^{n+2},Y^{j+1}\right\rangle =0\;\;(0\leq
j\leq n)$.
\end{proof}

\begin{remark}
Example \ref{exdegcaseIII} (below) illustrates a case of Proposition
\ref{degcaseIII} where $\mathcal{M}(2)$ admits no flat extension, so the
minimal representing measure is $(2n+2)$-atomic.
\end{remark}

\section{\label{proof}Proof of Theorem \ref{hyper}}

We now turn to the proof of Theorem \ref{hyper}, which we restate for the sake
of convenience. \ As in previous sections, it suffices to consider the cases
$yx=1$ and $yx=0$.

\begin{theorem}
\label{thmhyperbolic2}Let $\beta\equiv\beta^{(2n)}:\beta_{00},\beta_{01}%
,\beta_{10},...,\beta_{0,2n},...,\beta_{2n,0}$ be a family of real numbers,
$\beta_{00}>0$, and let $\mathcal{M}(n)$ be the associated moment matrix.
$\ $Assume that $\mathcal{M}(n)$ is positive, recursively generated, and
satisfies $YX=1$ (resp. $YX=0$). \ Then $\operatorname*{rank}\mathcal{M}%
(n)\leq2n+1$, and the following statements are equivalent.
\end{theorem}

\begin{enumerate}
\item[(i)] $\beta$\textit{\ admits a representing measure (necessarily
supported in }$yx=1$\textit{, resp. }$yx=0$\textit{).}

\item[(ii)] $\beta$\textit{\ admits a representing measure with convergent
moments up to degree }$2n+2$\textit{\ (necessarily supported in }%
$yx=1$\textit{, resp. }$yx=0$\textit{).}

\item[(iii)] $\beta$\textit{\ admits a representing measure }$\mu
$\textit{\ (necessarily supported in }$yx=1$\textit{, resp. }$yx=0$\textit{)
such that }$\operatorname*{card}\operatorname*{supp}\mu\leq
1+\operatorname*{rank}M(n)$\textit{. \ If }$\operatorname*{rank}M(n)\leq
2n$\textit{, then }$\mu$\textit{\ can be taken so that }$\operatorname*{card}%
\operatorname*{supp}\mu=\operatorname*{rank}M(n)$\textit{.}

\item[(iv)] $M(n)$\textit{\ admits a positive, recursively generated extension
}$M(n+1)$\textit{.}

\item[(v)] $M(n)$\textit{\ admits a positive, recursively generated extension
}$M(n+1)$\textit{, with }$\operatorname*{rank}M(n+1)\leq1+\operatorname*{rank}%
M(n)$\textit{, and }$M(n+1)$\textit{\ admits a flat extension }$M(n+2)$%
\textit{. \ If }$\operatorname*{rank}M(n)\leq2n$\textit{, then }%
$M(n)$\textit{\ admits a flat extension }$M(n+1)$\textit{.}

\item[(vi)] $\operatorname*{rank}\;M(n)\leq\operatorname*{card}\;V(\beta)$\textit{.}
\end{enumerate}

To establish Theorem \ref{thmhyperbolic2} we require the following result,
whose proof is an adaptation of the proof of Proposition \ref{degcaseI}.

\begin{proposition}
\label{propnew}Suppose $\mathcal{M}(n)\equiv\mathcal{M}(n)(\beta)$ is
positive, $YX=0$ in $\mathcal{C}_{\mathcal{M}(n)}$, and $\mathcal{M}(n)$ has a
positive, recursively generated extension $\mathcal{M}(n+1)$. \ Suppose also
that there exists $k$, $1<k\leq n$, such that $\mathcal{S}_{n}(k-1)$ is
linearly independent and $X^{k}\in$ lin.span $\mathcal{S}_{n}(k-1)$. \ Then
$\mathcal{M}(n)$ admits a flat extension $\mathcal{M}(n+1)$.
\end{proposition}

\begin{proof}
By hypothesis, we may write
\begin{equation}
X^{k}=a_{0}1+a_{1}X+b_{1}Y+...+a_{k-1}X^{k-1}+b_{k-1}Y^{k-1}\;\;(a_{i}%
,b_{i}\in\mathbb{R)}. \label{degneweq61}%
\end{equation}
Assume first that not all coefficients $b_{j}$ are zero, and let $m\leq k-1$
be the largest integer such that $b_{m}\neq0$. \ In $\mathcal{M}(n+1)$ we can
formally multiply (\ref{degneweq61}) by $Y$ to obtain
\[
0=YX^{k}=a_{0}Y+b_{1}Y^{2}+...+b_{m}Y^{m+1},
\]
from which it follows that $Y^{m+1}$ is a linear combination of columns
associated to powers of $y$ of lower degree, and, a fortiori, that the same is
true of $Y^{n}$. \ If we instead formally multiply (\ref{degneweq61}) by $X$,
we see that $X^{n}$ is a linear combination of columns associated to powers of
$x$ of lower degree. \ Thus, $\mathcal{M}(n)$ is flat. \ Since $\mathcal{M}%
(n+1)$ is a recursively generated extension of $\mathcal{M}(n)$, it must be a
flat extension of $\mathcal{M}(n)$, and the result follows in this case.

We can thus assume that all coefficients $b_{j}$ are zero, that is,
\begin{equation}
X^{k}=a_{0}1+a_{1}X+...+a_{k-1}X^{k-1}. \label{degneweq62}%
\end{equation}
If $a_{0}\neq0$, in $\mathcal{M}(n+1)$ we can\ formally multiply
(\ref{degneweq62}) by $Y$ to obtain $Y=0$, a contradiction. \ Thus,
(\ref{degneweq62}) does not involve the column $1$, just as in (\ref{degxk2}).
\ We may now continue exactly as in the proof of Proposition \ref{degcaseI},
since the part of that proof following (\ref{degxk2}) does not entail the
variety condition $\operatorname*{rank}\mathcal{M}(n)\leq\operatorname*{card}%
\mathcal{V}(\beta)$.
\end{proof}

\begin{proof}
[Proof of Theorem \textup{\ref{thmhyperbolic2}}]By \cite[Theorem 2.1]{tcmp2}
and the equivalence of the moment problems for $M(n)(\gamma)$ and
$\mathcal{M}(n)(\beta)$ \cite[Proposition 1.12]{tcmp6}, we can assume that the
columns $1,$ $X$ and $Y$ are linearly independent. \ Observe first that
$(iii)\Rightarrow(ii)\Rightarrow(i)$ trivially, that $(i)\Rightarrow(vi)$ by
\cite[(1.7)]{tcmp3}, and that $(iii)\Leftrightarrow(v)$ by (\ref{exp18}) and
(\ref{15}). \ Also, $(vi)\Rightarrow(iii)$ by Theorem \ref{firstmain}, so
$(i)$, $(ii)$, $(iii)$, $(v)$ and $(vi)$ are equivalent. \ Since
$(v)\Rightarrow(iv)$ is trivial, to complete the proof it suffices to
establish $(iv)\Rightarrow(v)$. \ 

Assume first that $\mathcal{M}(n)$ satisfies $YX=1$. \ The cases when
$\operatorname*{rank}\mathcal{M}(n)\leq2n$ correspond to the column dependence
relations that we considered in Cases I-III in the proof of Theorem
\ref{thmflat}, so we reconsider these dependence relations. \ First recall our
hypotheses: $\mathcal{M}(n)$ is positive, recursively generated, $YX=1$, and
$\mathcal{M}(n)$ admits a positive recursively generated extension
$\mathcal{M}(n+1)$. \ We need to show that $\mathcal{M}(n)$ admits a flat extension.

\textbf{Case I}. \ We have $\mathcal{S}_{n}(k-1)$ linearly independent and
$X^{k}=p_{k-1}(X)+q_{k-1}(Y)$ in $\mathcal{C}_{\mathcal{M}(n)}$, with $\deg
p_{k-1},$ $\deg q_{k-1}\leq k-1$. \ By the Extension Principle \cite{FiCM},
the same relation must hold in the column space of the positive extension
$\mathcal{M}(n+1)$. \ Since $\mathcal{M}(n+1)$ is recursively generated, we
must also have
\begin{equation}
X^{k-1}\equiv YX^{k}=Yp_{k-1}(X)+Yq_{k-1}(Y). \label{eq43}%
\end{equation}
\ Let $b_{k-1}$ be the coefficient of $y^{k-1}$ in $q_{k-1}$. \ If $b_{k-1}=0
$, then (\ref{eq43}) implies that $\mathcal{S}_{n}(k-1)${}is linearly
dependent, a contradiction. \ Thus, $b_{k-1}\neq0$, so (\ref{eq43}) implies
that $Y^{k}$ can be written as a linear combination of previous columns, and
therefore $\mathcal{M}(k)$ is flat; it now follows from \cite{FiCM} and
recursiveness that $\mathcal{M}(n)$ is flat, and thus admits a flat extension.

\textbf{Case II. \ }Suppose $\mathcal{S}:=\mathcal{S}_{n}(k-1)\bigcup
\{X^{k}\}$ is linearly independent and $Y^{k}\in$ lin.span $\mathcal{S}$, for
some $k<n$. \ The hypothesis about $\mathcal{M}(n+1)$ is superfluous, as we
showed in Proposition \ref{propcaseIII} that $\mathcal{M}(n)$ is flat.

\textbf{Case III}. \ Here $\mathcal{S}_{n}(n-1)\bigcup\{X^{n}\}$ is linearly
independent and
\begin{equation}
Y^{n}=a_{n}X^{n}+p_{n-1}(X)+q_{n-1}(Y) \label{eq7}%
\end{equation}
in $\mathcal{C}_{\mathcal{M}(n)}$, with $\deg p_{n-1},\deg q_{n-1}\leq n-1$.
\ The Extension Principle \cite{FiCM} shows that the same relation must hold
in the column space of the positive extension $\mathcal{M}(n+1)$. \ Since
$\mathcal{M}(n+1)$ is recursively generated, we must also have
\begin{equation}
Y^{n-1}\equiv Y^{n}X=a_{n}X^{n+1}+p_{n-1}(X)X+q_{n-1}(Y)X.\ \label{eq8}%
\end{equation}
If $a_{n}=0$, $\mathcal{S}_{n}(n-1)\bigcup\{X^{n}\}$ is linearly dependent, a
contradiction. \ Thus, $a_{n}\neq0$, which implies that $X^{n+1}$ can be
written in terms of columns of lower degree. \ Moreover, from (\ref{eq7}) we
obtain $Y^{n+1}=a_{n}X^{n-1}+Yp_{n-1}(X)+Yq_{n-1}(Y)$, so $Y^{n+1}$ is also a
linear combination of columns of lower degree. \ Finally, from $YX=1$ and the
recursiveness of $\mathcal{M}(n+1)$, we see that the intermediate columns
$Y^{i}X^{j}\;(i+j=n+1,$ with $i,j\geq1)$ are all identical to columns
corresponding to monomials of degree $n-1$. \ It follows that $\mathcal{M}%
(n+1)$ is flat, and is thus a flat extension of $\mathcal{M}(n)$.

\textbf{Case IV}. \ Here\ $\mathcal{S}_{n-1}(n)\bigcup\{X^{n}\}$we consider
the case when $\mathcal{S}_{n}(n)$ is linearly independent; the result follows
directly from Proposition \ref{propcaseV} (without using the given extension
$\mathcal{M}(n+1)$). \ This completes the proof for $YX=1$. \ 

We now assume that $\mathcal{M}(n)$ is positive, recursively generated, $YX=0
$ in $\mathcal{C}_{\mathcal{M}(n)}$, and $\mathcal{M}(n)$ admits a positive,
recursively generated extension $\mathcal{M}(n+1)$. \ We consider again the
various cases of column dependence relations that we examined in the proof of
Theorem \ref{secondflat}. \ If the first dependence relation in $\mathcal{S}%
_{n}(n)$ occurs at $X^{k}\;(2\leq k\leq n)$, then Proposition \ref{propnew}
implies that $\mathcal{M}(n)$ has a flat extension $\mathcal{M}(n+1)$. \ If
the first dependence relation occurs at $Y^{k}\;(2\leq k\leq n)$, then,
without recourse to the given extension $\mathcal{M}(n+1)$, Proposition
\ref{degcaseII} implies that $\mathcal{M}(n)$ admits a flat extension. \ In
the remaining case, $\operatorname*{rank}\mathcal{M}(n)=2n+1$, so the result
follows from Proposition \ref{degcaseIII} (again, without using the given
extension $\mathcal{M}(n+1)$).
\end{proof}

\section{\label{examples}Some Examples Illustrating Theorems \ref{thmflat} and
\ref{secondflat}}

Example \ref{ex16} illustrates Case III of Theorem \ref{thmflat}. \ We now
present examples corresponding to other cases of Theorems \ref{thmflat} and
\ref{secondflat}.

\begin{example}
\label{excasei}(Theorem \ref{thmflat}, Case I) We illustrate $\mathcal{M}(3)$
in which $YX=1$ and the first dependence relation in $\mathcal{S}_{3}(3)$
occurs at $X^{3}$. \ We define
\[
\mathcal{M}(3)\equiv\mathcal{M}(3)(\beta):=\left(
\begin{array}
[c]{cccccccccc}%
1 & 0 & 0 & a & 1 & a & 0 & 0 & 0 & 0\\
0 & a & 1 & 0 & 0 & 0 & 2a^{2} & a & 1 & a\\
0 & 1 & a & 0 & 0 & 0 & a & 1 & a & 2a^{2}\\
a & 0 & 0 & 2a^{2} & a & 1 & e & 0 & 0 & 0\\
1 & 0 & 0 & a & 1 & a & 0 & 0 & 0 & 0\\
a & 0 & 0 & 1 & a & 2a^{2} & 0 & 0 & 0 & f\\
0 & 2a^{2} & a & e & 0 & 0 & g & 2a^{2} & a & 1\\
0 & a & 1 & 0 & 0 & 0 & 2a^{2} & a & 1 & a\\
0 & 1 & a & 0 & 0 & 0 & a & 1 & a & 2a^{2}\\
0 & a & 2a^{2} & 0 & 0 & f & a & a & 2a^{2} & h
\end{array}
\right)  ,
\]
where $a>1$. \ Clearly $YX=1$, and a calculation shows that
$\operatorname*{rank}\mathcal{M}(2)=5$. \ Thus, the set $\mathcal{B}%
:=\{1,X,Y,X^{2},Y^{2}\}$ is linearly independent, and
\begin{equation}
\lbrack X^{3}]_{\mathcal{B}}=-ae[1]_{\mathcal{B}}+\frac{a(2a^{2}-1)}{a^{2}%
-1}[X]_{\mathcal{B}}-\frac{a^{2}}{a^{2}-1}[Y]_{\mathcal{B}}+\frac{a^{2}%
e}{2a^{2}-1}[X^{2}]_{\mathcal{B}}+\frac{e(a^{2}-1)}{2a^{2}-1}[Y^{2}%
]_{\mathcal{B}}. \label{eqx}%
\end{equation}
We seek to impose the same relation in $\mathcal{C}_{\mathcal{M}(3)}$, which
entails
\[
g=\frac{a(2a^{2}-1)}{a^{2}-1}2a^{2}-\frac{a^{2}}{a^{2}-1}a+\frac{a^{2}%
e}{2a^{2}-1}e=\allowbreak a^{2}\frac{8a^{5}-10a^{3}+3a+a^{2}e^{2}-e^{2}%
}{\left(  a^{2}-1\right)  \left(  2a^{2}-1\right)  }%
\]
and
\[
\frac{a(2a^{2}-1)}{a^{2}-1}a-\frac{a^{2}}{a^{2}-1}2a^{2}+\frac{e(a^{2}%
-1)}{2a^{2}-1}f=a,
\]
so that
\[
f=a\frac{(a^{2}+a-1)(2a^{2}-1)}{e\left(  a^{2}-1\right)  ^{2}}.
\]
To compute $\mathcal{V}(\mathcal{M}(3))$, note that the column relation
(\ref{eqx}), together with $YX=1$, gives rise to the equation
\[
x^{3}+ae-\frac{a(2a^{2}-1)}{a^{2}-1}x+\frac{a^{2}}{a^{2}-1}\frac{1}%
{x}-\frac{a^{2}e}{2a^{2}-1}x^{2}-\frac{e(a^{2}-1)}{2a^{2}-1}\frac{1}{x^{2}%
}=0,
\]
or, equivalently,%
\begin{gather}
(2a^{2}-1)(a^{2}-1)x^{5}-a^{2}e(a^{2}-1)x^{4}-a(2a^{2}-1)^{2}x^{3}%
\label{equa}\\
+ae(a^{2}-1)(2a^{2}-1)x^{2}+a^{2}(2a^{2}-1)x-e(a^{2}-1)^{2}=0.\nonumber
\end{gather}
Thus, $\mathcal{V}(\mathcal{M}(3))$ can have at most five points, so if
$\operatorname*{rank}\mathcal{M}(3)\leq\operatorname*{card}\mathcal{V}%
(\mathcal{M}(3))$, then column $Y^{3}$ must be written in terms of the columns
in $\mathcal{B}$ (in particular, $h$ is fully determined), and $\mathcal{M}%
(3)$ is therefore flat, in accordance with Proposition \ref{propcaseI}. \ For
a specific numerical example, let $a:=2$ and $e:=1$, so that $f=\frac{70}{9}$,
$g=\frac{740}{21}$, $h=\frac{2190400}{3087}$ and $Y^{3}=-\frac{1480}{21}%
\cdot1-\frac{4}{3}X+\frac{14}{3}Y+\frac{740}{49}X^{2}+\frac{2960}{147}Y^{2}$.
\ Then (\ref{equa}) becomes $21x^{5}-12x^{4}-98x^{3}+42x^{2}+28x-9=0$, which
has five real roots, as follows: $x_{1}\cong-2.0292$, $x_{2}\cong-0.521229$,
$x_{3}\cong0.282006$, $x_{4}\cong0.658788$ and $x_{5}\cong2.18106$. \ To
calculate the densities we use Theorem \ref{thm17}; here $V$ is the $5\times5$
matrix whose entry in row $k$, column $\ell$ is $y_{\ell}^{i_{k}}x_{\ell
}^{j_{k}}\;(1\leq k,\ell\leq r)$, where $(i_{1},j_{1})=(0,0)$, $(i_{2}%
,j_{2})=(0,1)$, $(i_{3},j_{3})=(1,0)$, $(i_{4},j_{4})=(0,2)$ and $(i_{5}%
,j_{5})=(2,0)$. $\ $We then solve the equation $V\rho^{t}=(\beta_{i_{1},j_{1}%
},...,\beta_{i_{5},j_{5}})^{t}$, where $\rho\equiv(\rho_{1},...,\rho_{5})$ and
$(\beta_{i_{1},j_{1}},...,\beta_{i_{5},j_{5}})=(1,0,0,2,2)$. \ Thus, $\rho
_{1}\cong0.228429$, $\rho_{2}\cong0.263185$, $\rho_{3}\cong0.0174322$,
$\rho_{4}\cong0.31204$, $\rho_{5}\cong0.178914$. \ \textbf{\qed                }
\end{example}

\begin{example}
(Theorem \ref{thmflat}, Cases III and IV) \ Consider the moment matrix
\[
\mathcal{M}(3)\equiv\mathcal{M}(3)(\beta):=\left(
\begin{array}
[c]{cccccccccc}%
1 & 0 & 0 & 1 & 1 & 2 & 0 & 0 & 0 & 0\\
0 & 1 & 1 & 0 & 0 & 0 & 3 & 1 & 1 & 2\\
0 & 1 & 2 & 0 & 0 & 0 & 1 & 1 & 2 & 5\\
1 & 0 & 0 & 3 & 1 & 1 & 0 & 0 & 0 & 0\\
1 & 0 & 0 & 1 & 1 & 2 & 0 & 0 & 0 & 0\\
2 & 0 & 0 & 1 & 2 & 5 & 0 & 0 & 0 & 3\\
0 & 3 & 1 & 0 & 0 & 0 & 14 & 3 & 1 & 1\\
0 & 1 & 1 & 0 & 0 & 0 & 3 & 1 & 1 & 2\\
0 & 1 & 2 & 0 & 0 & 0 & 1 & 1 & 2 & 5\\
0 & 2 & 5 & 0 & 0 & 3 & 1 & 2 & 5 & 33
\end{array}
\right)  .
\]
It is easy to see that $\mathcal{M}(3)\geq0$, $YX=1$ in $\mathcal{C}%
_{\mathcal{M}(3)}$, $\operatorname*{rank}\mathcal{M}(3)=7$, and that
$\mathcal{S}_{3}(3)$ is a basis for $\mathcal{C}_{\mathcal{M}(3)}$. \ The
block $B(4)$ for a recursively generated extension is determined by the choice
of $\beta_{07}=t$ and $\beta_{70}=u$. \ Following the proof of Proposition
\ref{propcaseV}, we have Ran $B(4)\subseteq$ Ran $\mathcal{M}(3)$, and a
calculation of $\mathcal{M}\equiv\lbrack\mathcal{M}(3);B(4)]$ shows that
$\mathcal{M}$ has the form of a moment matrix $\mathcal{M}(4)$ if and only if
$(t+3)(u-150)=1$. \ For example, with $t=-2$ and $u=151$, we find the
following column relations in $\mathcal{M}(4)$ that are not determined
recursively from relations in $\mathcal{M}(3)$:
\[
X^{4}=-20\cdot1+4X+Y+9X^{2}+7Y^{2}-X^{3}-Y^{3}%
\]
and
\[
Y^{4}=4\cdot1+9X-20Y-X^{2}+Y^{2}-X^{3}+7Y^{3}.\
\]
Together with $YX=1$, these relations show that $\operatorname*{card}%
\mathcal{V}(\mathcal{M}(4))=\operatorname*{rank}\mathcal{M}(4)=7$; the
$x$-coordinates of the points in $\{(x_{i},\frac{1}{x_{i}})\}_{i=1}^{7}%
\equiv\mathcal{V}(\mathcal{M}(4))$ are as follows: $x_{1}\cong-2.82238$,
$x_{2}\cong-1.87650$, $x_{3}\cong-0.64947$, $x_{4}\cong0.14873$, $x_{5}%
\cong0.66898$, $x_{6}\cong1.32445$ and $x_{7}\cong2.20619$. \ Since
$\mathcal{S}_{3}(3)$ is a basis for $\mathcal{C}_{\mathcal{M}(3)}$, Theorem
\ref{thm17} now implies that the densities $\{\rho_{i}\}_{i=1}^{7}$ of the
unique representing measure $\mu=\sum_{i=1}^{7}\rho_{i}\delta_{(x_{i}%
,\frac{1}{x_{i}})}$ for $\mathcal{M}(4)$ are as follows: $\rho_{1}%
\cong0.00754$, $\rho_{2}\cong0.07579$, $\rho_{3}\cong0.41491$, $\rho_{4}%
\cong0.00025$, $\rho_{5}\cong0.43154$, $\rho_{6}\cong0.01146$ and $\rho
_{7}\cong0.05851$.
\end{example}

We can illustrate the second part of Case IV and also Case III if we use
$t=-3$ and $u=150$, so that $M$ is not a moment matrix. \ In this case, the $C
$ block of $M$ satisfies $C_{11}=79$. \ Following the proof of Proposition
\ref{propcaseV}, we redefine $C_{11}:=80$, $C_{15}:=1$, and we compute
$C_{55}\;(=\beta_{80})=1036$, so as to define a positive rank $8$ moment
matrix $M(4)$ in which the first dependence relation in $S_{4}(4)$ is of the
form $Y^{4}=-16\cdot1+18X-21Y+8X^{2}+8Y^{2}-3X^{3}+6Y^{3}-X^{4}$ (i.e., $M(4)$
is in Case III). \ This relation and $YX=1$ determine $\mathcal{V}%
\equiv\mathcal{V}(\mathcal{M}(4))$, and a calculation shows that
$\operatorname*{card}\mathcal{V}=\operatorname*{rank}\mathcal{M}(4)=8$; the
$x$-coordinates of the points $\{(x_{i},\frac{1}{x_{i}})\}_{i=1}^{8}%
\equiv\mathcal{V}$ are as follows: $x_{1}\cong-3.63582$, $x_{2}\cong-2.02578$,
$x_{3}\cong-0.840968$, $x_{4}\cong-0.644996$, $x_{5}\cong0.149637$,
$x_{6}\cong0.670287$, $x_{7}\cong1.14158$ and $x_{8}\cong2.18606$.\ \ Since
the unique representing measure $\mu$ for the flat extension $M(5)$ of $M(4)$
(guaranteed by Proposition \ref{propcaseIV}) satisfies $\operatorname*{card}%
\operatorname*{supp}\mu=\operatorname*{rank}M(4)=8$, it follows that
$\operatorname*{supp}\mu=\mathcal{V}$ and that the densities of $\mu$
(computed as in Theorem \ref{thm17}) are $\rho_{1}\cong0.00087$, $\rho
_{2}\cong0.07462$, $\rho_{3}\cong0.02883$, $\rho_{4}\cong0.39396$, $\rho
_{5}\cong0.00025$, $\rho_{6}\cong0.43357$, $\rho_{7}\cong0.00612$ and
$\rho_{8}\cong0.06177 $. \textit{\qed                }

We next turn to examples illustrating Theorem \ref{secondflat}.

\begin{example}
(cf. Proposition \ref{degcaseI}) \ We illustrate $\mathcal{M}(2)$ with $YX=0$,
where the first dependence relation in $\mathcal{S}_{2}(2)$ occurs at $X^{2}$.
\ Consider the moment matrix
\[
\mathcal{M}(2):=\left(
\begin{array}
[c]{cccccc}%
1 & a & 0 & b & 0 & d\\
a & b & 0 & e & 0 & 0\\
0 & 0 & d & 0 & 0 & g\\
b & e & 0 & f & 0 & 0\\
0 & 0 & 0 & 0 & 0 & 0\\
d & 0 & g & 0 & 0 & h
\end{array}
\right)  .
\]
A straightforward calculation reveals that with $b>a^{2}$, $d>0$,
$e=\frac{b^{2}}{a}$, $f=\frac{b^{3}}{a^{2}}$ and $h>\frac{bd^{3}-a^{2}%
g^{2}+bg^{2}}{d(b-a^{2})}$, $\mathcal{M}(2)\geq0$, $\operatorname*{rank}%
\mathcal{M}(2)=4$, $YX=0 $ in $\mathcal{C}_{\mathcal{M}(2)}$, and
$X^{2}=\frac{b}{a}X$. \ The variety $\mathcal{V}(\mathcal{M}(2))$ associated
to this matrix consists of the intersection of $yx=0$ with $x(x-\frac{b}%
{a})=0$, so $\operatorname*{card}\mathcal{V}(\mathcal{M}(2))=\infty$. \ Now,
in any recursively generated extension $\mathcal{M}(3)$ we must have
$X^{3}=\frac{b}{a}X^{2}=\frac{b^{2}}{a^{2}}X$. \ Also, the form of a block
$B(3)$ is
\[
B(3)\equiv\left(
\begin{array}
[c]{cccc}%
e & 0 & 0 & g\\
f & 0 & 0 & 0\\
0 & 0 & 0 & h\\
p & 0 & 0 & 0\\
0 & 0 & 0 & 0\\
0 & 0 & 0 & q
\end{array}
\right)  ,
\]
where $p:=\frac{b^{4}}{a^{3}}$ and $q$ is free. \ Following the proof of
Proposition \ref{degcaseI}, the existence of a flat extension $\mathcal{M}(3)
$ depends on the verification of (\ref{eqnew8}), that is, $[X^{3}%
]_{\{1,2,3,6\}}\mathcal{M}(2)_{\{1,2,3,6\}}^{-1}[Y^{3}]_{\{1,2,3,6\}}=0 $. \ A
straightforward calculation shows that this is indeed the case, so
$\mathcal{M}(2)$ admits a flat extension. \textbf{\qed                }
\end{example}

\begin{example}
\label{exdegcaseIII}(cf. Proposition \ref{degcaseIII} and \cite[Example
5.6]{tcmp6})\ Consider the moment matrix
\[
\mathcal{M}(2):=\left(
\begin{array}
[c]{cccccc}%
1 & 1 & 1 & 2 & 0 & 3\\
1 & 2 & 0 & 4 & 0 & 0\\
1 & 0 & 3 & 0 & 0 & 9\\
2 & 4 & 0 & 9 & 0 & 0\\
0 & 0 & 0 & 0 & 0 & 0\\
3 & 0 & 9 & 0 & 0 & 28
\end{array}
\right)
\]
Observe that $\operatorname*{rank}\mathcal{M}(2)=5$, and that $YX=0$ in
$\mathcal{C}_{\mathcal{M}(2)}$, so that $\mathcal{S}_{2}(2)$ is a basis for
$\mathcal{C}_{\mathcal{M}(2)}$. \ A block $B(3)$ for a recursively generated
extension is completely determined by the choice of $\beta_{05}=p$ and
$\beta_{50}=q$. \ With these choices, the $C$ block in the flat extension
$[\mathcal{M}(2);B(3)]$ has entries $C_{11}=(p-18)^{2}+42$, $C_{14}=C_{41}=1$,
and $C_{44}=(q-84)^{2}+262$. \ It is then clear that $[\mathcal{M}(2);B(3)]$
is not a moment matrix, so $\beta^{(4)}$ has no $5$-atomic representing
measure. \ To construct a $6$-atomic representing measure, we modify the two
key entries, $C_{11}$ and $C_{14}$, as dictated by the proof of Proposition
\ref{degcaseIII}. \ By taking $p=18$, $q=84$, $u=43$, it is easy to see that
$\operatorname*{rank}\mathcal{M}(3)=1+\operatorname*{rank}\mathcal{M}%
(2)\;(=6)$ precisely when $v=263$. \ As in the proof of Proposition
\ref{degcaseIII}, we claim that $\mathcal{M}(3)$ admits a flat extension
$\mathcal{M}(4)$. \ We first exhibit a unique block $B(4)$, subordinate to
$YX=0$, such that Ran $B(4)\subseteq$ Ran $\mathcal{M}(3)$. \ In any such
block, $YX^{3}=Y^{2}X^{2}=Y^{3}X=0$. \ In $Y^{4}$, all the entries are
determined from $\mathcal{M}(3)$ and $YX=0$, except $s:=\left\langle
Y^{4},Y^{3}\right\rangle $. \ Now, in $\mathcal{C}_{\mathcal{M}(3)}$ we have a
dependence relation
\[
Y^{3}=-5\cdot1+7X+11Y-X^{3};
\]
since a flat extension must necessarily be recursively generated, in
$\mathcal{C}_{\mathcal{M}(4)}$ we must have
\[
Y^{4}=-5Y+11Y^{2},
\]
whence $s=\left\langle -5Y+11Y^{2},Y^{3}\right\rangle =-5\beta_{4,0}%
+11\beta_{5,0}=784$. \ Also, in $X^{4}$, all values are determined except
$r:=\left\langle X^{4},X^{3}\right\rangle $; following Proposition
\ref{degcaseIII}, we claim that there is a unique value for $r$ such that
$X^{4}\in$ Ran $\mathcal{M}(3)$. \ For, relative to the linearly independent
set of columns $\{1,X,Y,X^{2},Y^{2},X^{3}\}$ in $\mathcal{C}_{\mathcal{M}(3)}%
$,
\[
X^{4}=(3r-243)\cdot1+(481-6r)X+(81-r)Y+7X^{2}+(r-81)X^{3}\text{,}%
\]
and since $\left\langle X^{4},Y^{3}\right\rangle $ must be zero, we obtain
\begin{align*}
0  &  =\left\langle (3r-243)\cdot1+(481-6r)X+(81-r)Y+7X^{2}+(r-81)X^{3}%
,Y^{3}\right\rangle \\
&  =(3r-243)\beta_{3,0}+(81-r)\beta_{4,0}=(3r-243)\cdot9+(81-r)\cdot28=81-r,
\end{align*}
whence $r=81$. \ The unique $6$-atomic representing measure associated with
$\mathcal{M}(3)$ is then given by $(x_{1},y_{1})\cong(2.16601,0)$,
$(x_{2},y_{2})\cong(0.782816,0)$, $(x_{3},y_{3})\cong$\newline $(-2.94883,0)$,
$(x_{4},y_{4})\cong(0,0.463604)$, $(x_{5},y_{5})\cong(0,3.06043)$,
$(x_{6},y_{6})\cong$\newline $(0,-3.52404)$, $\rho_{1}\cong0.393081$,
$\rho_{2}\cong0.203329$, $\rho_{3}\cong0.00359018$, $\rho_{4}\cong0.0821253$,
$\rho_{5}\cong0.316218$, $\rho_{6}\cong0.00165656$. \textbf{\qed              }
\end{example}

\section{\label{application}An Application to the Full Moment Problem}

We conclude with a new proof of Theorem \ref{thm13}. \ We require the
following preliminary result, which was given in \cite{FiaNew} for complex
moment matrices, but also holds for $\mathcal{M}\left(  \infty\right)
(\beta)$.

\begin{lemma}
\label{lemfia}(\cite[Proposition 4.2]{FiaNew}) \ Let $\beta\equiv
\beta^{(\infty)}$ be a full sequence such that $\mathcal{M}\left(
\infty\right)  \geq0$. \ Then $\mathcal{M}(n)$ is positive and recursively
generated for each $n\geq1$.
\end{lemma}

\begin{theorem}
(cf. \cite{Sto1}) \ Let $P\in\mathbb{R}[x,y]$ with $\deg P\leq2$. \ The
sequence $\beta\equiv\beta^{(\infty)}$ has a representing measure supported in
$P(x,y)=0$ if and only if $\mathcal{M}\equiv\mathcal{M}(\infty)(\beta)\geq0$
and $P(X,Y)=0$ in $\mathcal{C}_{\mathcal{M}}$.
\end{theorem}

\begin{proof}
Let $\mu$ be a representing measure for $\beta$ supported in $\mathcal{Z}(P)$.
\ For each $n\geq2$, $\mu$ is a representing measure for $\beta^{(2n)}$, so
(\ref{14}) and (\ref{15}) imply that $\mathcal{M}(n)(\beta)\geq0$ and
$P(X,Y)=0$ in $\mathcal{C}_{\mathcal{M}(n)}$; thus $\mathcal{M}\geq0$ and
$P(X,Y)=0$ in $\mathcal{C}_{\mathcal{M}}$.

For the converse, since $\mathcal{M}\geq0$, Lemma \ref{lemfia} implies that
for $n\geq2$, $\mathcal{M}(n+1)$ is positive and recursively generated, so
$\mathcal{M}(n)$ has a positive, recursively generated extension. \ Since
$P(X,Y)=0$ in $\mathcal{C}_{\mathcal{M}}$, the same is true in $\mathcal{C}%
_{\mathcal{M}(n)}$. \ Suppose now that $P(x,y)=0$ is a nondegenerate
hyperbola. \ Let $\Phi$ be an injective degree-one mapping of the plane onto
itself such that $\Phi(\mathcal{Z}(P))=\{(x,y):yx=1\}$ (cf. Section
\ref{Int}), and let $\tilde{\beta}^{(2n)}$ be the sequence corresponding to
$\beta^{(2n)}$ via Proposition \ref{prop17}. \ Proposition \ref{prop17}(ix)
implies that $\mathcal{\tilde{M}}(n)\equiv\mathcal{M}(n)(\tilde{\beta})$ has a
positive, recursively generated extension $\mathcal{\tilde{M}}(n+1)$, and
Proposition \ref{prop17}(vi) implies that $YX=1$ in $\mathcal{C}%
_{\mathcal{\tilde{M}}(n)}$. \ Theorem \ref{hyper}(iv) now implies that
$\tilde{\beta}^{(2n)}$ has a representing measure supported in $yx=1$, so
Proposition \ref{prop17}(v) implies that $\beta^{(2n)}$ has a representing
measure supported in $P(x,y)=0$; the existence of a representing measure for
$\beta$ supported in $P(x,y)=0$ now follows from Theorem \ref{thm14}.

The cases of other conics are handled similarly. \ For degenerate hyperbolas,
we use the case of Theorem \ref{firstmain} for $yx=0$. \ For parabolas and
ellipses we use \cite[Theorem 1.4]{tcmp7} and \cite[Theorem 3.5]{tcmp5},
respectively. \ The cases of lines is treated directly in \cite[Theorem
2.1]{tcmp2}.
\end{proof}

\begin{remark}
(i) By analogy with properties $(A)$ and $(A_{n}^{\prime})$ (cf. Section
\ref{Int}), consider also the following possible properties for a polynomial
$P\in\mathbb{R}[x,y]$:\
\begin{gather}
\beta\equiv\beta^{(2n)}\text{ has a representing measure supported in
}\mathcal{Z}(P)\text{ if and only if }\nonumber\\
\mathcal{M}(n)(\beta)\text{ is positive semi-definite, recursively generated
}\tag{$A_n$}\\
\text{and }P(X,Y)=0\text{ in }\mathcal{C}_{\mathcal{M}(n)}.\nonumber
\end{gather}%
\begin{gather}
\beta\equiv\beta^{(\infty)}\text{ has a representing measure supported in
}\mathcal{Z}(P)\text{ if and only if }\nonumber\\
\mathcal{M}(\infty)(\beta)\text{ is positive semi-definite, }P(X,Y)=0\text{ in
}\mathcal{C}_{\mathcal{M}(\infty)}\text{, }\tag{$A^\prime$}\\
\text{and }\operatorname*{rank}\mathcal{M}(\infty)\leq\operatorname*{card}%
\mathcal{V}(\mathcal{M}(\infty)).\nonumber
\end{gather}
Theorem \ref{thm14} readily implies the following implications concerning
possible properties enjoyed by a given polynomial $P$:
\[%
\begin{array}
[c]{ccc}%
(A_{n})\text{ holds for all }n\geq\deg P & \Rightarrow & (A_{n}^{\prime
})\text{ holds for all }n\geq\deg P\\
\Downarrow &  & \Downarrow\\
(A)\text{ holds} & \Rightarrow & (A^{\prime})\text{ holds.}%
\end{array}
\]
First-degree polynomials $P$ satisfy $(A_{n})$ for all $n\geq\deg P$, but our
results show that second-degree polynomials do not. \ However, second-degree
polynomials satisfy $(A_{n}^{\prime})$ for all $n\geq\deg P$, and also $(A)$.
\ Stochel \cite{Sto1} has identified cubics which fail to satisfy $(A)$, but
we know of no example of a polynomial $P$ that fails to satisfy $(A^{\prime}%
)$, or $(A_{n}^{\prime})$ for some $n\geq\deg P$. \ \newline (ii) The full
moment problem on compact semi-algebraic sets in $\mathbb{R}^{n}$ was solved
by K. Schm\"{u}dgen \cite{Sch}. \ Recently, the analysis of the semi-algebraic
case was extended to non-compact sets by V. Powers and C. Scheiderer
\cite{PoS} (cf. \cite{KuM}, \cite{Sche}, \cite{SchNew}).
\end{remark}

\end{document}